\documentclass[11pt]{article}
\usepackage{amssymb}
\usepackage{amsmath}
\usepackage{array}
\usepackage[a4paper]{geometry}
\usepackage{lscape}

\advance \topmargin by -\headheight
\advance \topmargin by -\headsep
\evensidemargin \oddsidemargin
\input tcilatex
\newtheorem{thm}{Theorem}

\newtheorem{assu}{Assumption}
\newtheorem{lem}{Lemma}
\newtheorem{prop}{Proposition}

\begin{document}

\author{Jushan Bai \thanks{%
Department of Economics, New York University, New York, NY 10003,USA, and
School of Economics and Management, Tsinghua University. Email:
Jushan.Bai@nyu.edu.} \and Chihwa Kao\thanks{
Center for Policy Research and Department of Economics, Syracuse University,
Syracuse, NY\ 13244-1020, USA. Email: cdkao@maxwell.syr.edu.} \and Serena Ng
\thanks{
Department of Economics, Columbia University, 440 W. 118 St. New York,
NY 10027, USA,
Email: serena.Ng@columbia.edu. \newline
We thank Yixiao Sun, Joon Park, Yoosoon Chang, and Joakim Westerlund for
helpful discussions. We would also like to thank seminar and conference
participants at NYU, Texas A\&M, CESG 2006, Academia Sinica, National Taiwan
University, Greater New York Area Econometrics Colloquium at Yale, and CIREQ
Time Series Conference in Montreal for helpful comments. Bai and Ng
gratefully acknowledge financial support from the NSF (grants SES-0551275
and SES-0549978). } }
\title{\textbf{PANEL COINTEGRATION WITH GLOBAL STOCHASTIC TRENDS}}
\date{This version: \today }
\maketitle

\begin{abstract}
This paper studies estimation of panel cointegration models with
cross-sectional dependence generated by unobserved global stochastic trends.
The standard least squares estimator is, in general, inconsistent owing to
the spuriousness induced by the unobservable I(1) trends. We propose two
iterative procedures that jointly estimate the slope parameters and the
stochastic trends. The resulting estimators are referred to respectively as
CupBC (continuously-updated and bias-corrected) and the CupFM
(continuously-updated and fully-modified) estimators. We establish their
consistency and derive their limiting distributions. Both are asymptotically
unbiased and asymptotically normal and permit inference to be conducted
using standard test statistics. The estimators are also valid when there are
mixed stationary and non-stationary factors, as well as when the factors are
all stationary.

\vspace*{.5in} \noindent \textit{JEL Classification}: C13; C33 \bigskip

\end{abstract}

\thispagestyle{empty} \setcounter{page}{0} \newpage \baselineskip=18.0pt

\section{Introduction}

This paper is concerned with estimating panel cointegration models using a
large panel of data. Our focus is on estimating the slope parameters of the
non-stationary regressors when the cross sections share common sources of
non-stationary variation in the form of global stochastic trends. The
standard least squares estimator is either inconsistent or has a slower
convergence rate. We provide a framework for estimation and inference. We
propose two iterative procedures that estimate the latent common trends
(hereafter factors) and the slope parameters jointly. The estimators are $%
\sqrt{n}T$ consistent and asymptotically mixed normal. As such,
inference can be made using standard $t$ and Wald tests. The
estimators are also valid when some or all of the common factors are
stationary, and when some of the observed regressors are stationary.

Panel data have long been used to study and test economic
hypotheses. Panel data bring in information from two dimensions
to permit analysis that would otherwise be inefficient, if not impossible, with
time series or cross-sectional data alone. A new development in recent years
is the use of `large dimensional panels', meaning that the sample size in
the time series ($T$) and the cross-section $(n)$ dimensions are both large.
This is in contrast to traditional panels in which we have data of many
units over a short time span, or of a few variables over a long horizon.
Many researchers have come up with new ideas to exploit the rich information
in large panels.\footnote{%
See, for example, Baltagi (2005), Hsiao (2003), Pesaran and Smith (1995),
Kao (1999), and Moon and Phillips (2000, 2004) in the context of testing the
unit root hypothesis using panel data. Stock and Watson (2002) suggest
diffusion-index forecasting, while Bernanke and Boivin (2003) suggest new
formulations of vector autoregressions to exploit the information in large
panels.} However, large panels also raise econometric issues of their own.
In this analysis, we tackle two of these issues: the data $(y_{it},x_{it})$
are non stationary, and the structural errors $e_{it}=y_{it}-x_{it}^{\prime
}\beta $ are neither iid across $i$ nor over $t$. Instead, they are
cross-sectionally dependent and strongly persistent and possibly
non-stationary. In addition, $e_{it}$ are also correlated with the
explanatory variables $x_{it}$. These problems are dealt with by putting a
factor structure on $e_{it}$ and modelling the factor process explicitly.

The presence of common sources of non-stationarity leads naturally to the
concept of cointegration. In a small panel made up of individually I(1) (or
unit root) processes $y_t$ and $x_t$, where small means that the dimension
of $y_t$ plus the dimension of $x_t$ is treated as fixed in asymptotic
analysis, cointegration as defined in Engle and Granger (1986) means that
there exists a cointegrating vector, $(1 \; -\beta^\prime)$, such that the
linear combinations $y_t-x_t^\prime\beta$ are stationary, or are an I(0)
processes. In a panel data model specified by $y_{it}=x_{it}^\prime \beta +
e_{it}$ where $y_{it}$ and $x_{it}$ are I(1) processes, and that $e_{it}$
are iid across $i$, cointegration is said to hold if $e_{it}$ are `jointly'
I(0), or in other words, $(1,\; -\beta)$ is the common cointegrating vector
between $y_{it}$ and $x_{it}$ for \emph{all} $n$ units. A large literature
 panel cointegration already exists\footnote{%
See, for example, Phillips and Moon (1999) and Kao (1999). Recent surveys
can be found in Baltagi and Kao (2000) and Breitung and Pesaran (2005).} for
modelling panel cointegration when $e_{it}$ is cross-sectionally independent.
A serious drawback of  panel cointegration models with cross-section
independence is that there is  no role
for  common shocks, which, in theory should be the underlying source of
comovement in the cross-section units. Failure to account for common shocks
can potentially invalidate estimation and inference of $\beta $.\footnote{%
Andrews (2005) showed that cross-section dependence induced by common shocks
can yield inconsistent estimates. Andrews' argument is made in the context
of a single cross section and for stationary regressors and errors. For a
single cross section, not much can be done about common shocks. But for
panel data, we can explore the common shocks to yield consistent procedures.}
In view of this, more recent work has allowed for cross-sectional dependence
of $e_{it}$ when testing for the null hypothesis of panel cointegration.%
\footnote{%
See, for example, Phillips and Sul (2003), Gengenbach et al. (2005b), and
Westerlund (2006).} There is also a growing literature on panel unit root
tests with cross-sectional dependence.\footnote{%
For example, Chang (2002,2004), Choi (2006), Moon and Perron (2004),
Breitung and Das (2005), Gengenbach et al. (2005a), and Westerlund and
Edgerton (2006). Breitung and Pesaran (2005) provide additional references
in their survey.}


In this paper, we consider estimation and inference of
parameters in a panel model with cross-sectional dependence in the form of
common stochastic trends. The framework we adopt is that $e_{it}$ has
a common  component and a
stationary idiosyncratic component. That is, $e_{it}=\lambda _{i}^{\prime
}F_{t}+u_{it}$, so that panel cointegration holds when $u_{it}=y_{it}-\beta
x_{it}-\lambda _{i}^{\prime }F_{t}$ is jointly stationary.  A
regression of $y_{it}$  on $x_{it}$ will give a consistent estimator for $%
\beta $ when $F_{t}$ is I(0). We focus on estimation and inference
about $\beta $ when  $F_{t}$ is non-stationary. Empirical studies
suggest the relevance of such a setup. Holly, Pesaran and Yamagata
(2006) analyzed the  relationship  between real housing
prices and real income at state level, allowing for unobservable common
factors. They found evidence of cointegration after controlling
for common factors and additional spatial correlations. Some economic models
lead naturally to this set up. Consider a panel of industry data on
output, factor inputs such as  capital, and labor. Neoclassical
production function suggests that log output $y_{it}$ is linear in
log factor inputs $x_{it}$ and log productivity $e_{it}$.
Decomposing the latent $e_{it}$ into  the industry wide component
$F_t$ and an industry specific component $u_{it}$ and assuming that
$F_t$ is the source of non-stationarity leads to a model with latent
common trends.   In such a case, a regression of $y_{it}$ on
$x_{it}$ is spurious since $e_{it}$ is not only cross-sectionally
correlated, but also non-stationary.




 We deal with the problem by treating the common I(1)
variables as parameters. These are estimated jointly with $\beta $ using an
iterated procedure. The procedure is shown to yield a consistent estimator
of $\beta $, but the estimator is asymptotically biased. We then construct
two estimators to account for the bias arising from endogeneity and serial
correlation so as to re-center the limiting distribution around zero. The
first, denoted CupBC, estimates the asymptotic bias directly. The second,
denoted CupFM, modifies the data so that the limiting distribution does not
depend on nuisance parameters. Both are `continuously updated' (Cup)
procedures and require iteration till convergence. The estimators are $\sqrt{%
n}T$ consistent for the common slope coefficient vector, $\beta $. The
estimators enable use of standard test statistics such as $t$, $F$, and $%
\chi ^{2}$ for inference. The estimators are robust to mixed I(1)/I(0)
factors, as well as mixed I(1)/I(0) regressors. Thus, our approach is an
alternative to the solution proposed in Bai and Kao (2006) for stationary
factors. As we argue below, the Cup estimators have some advantages that
make an analysis of their properties interesting in its own right.

The rest of the paper is organized as follows. Section 2 describes the basic
model of panel cointegration with unobservable common stochastic trends.
Section 3 develops the asymptotic theory for the continuously-updated and
fully-modified estimators. Section 4 examines issues related to incidental
trends, mixed I(0)/I(1) regressors and mixed I(0)/I(1) common shocks, and
issues of testing cross-sectional independence. Section 5 presents Monte
Carlo results to illustrate the finite sample properties of the proposed
estimators. Section 6 provides a brief conclusion. The appendix contains the
technical materials.


\section{The Model}

Consider the model
\begin{equation*}
y_{it}=x_{it}^{\prime }\beta +e_{it}
\end{equation*}%
where for $i=1,...,n,$ $t=1,...,T,$ $y_{it}$ is a scalar,
\begin{equation}
x_{it}=x_{it-1}+\varepsilon _{it}.  \label{xit}
\end{equation}%
is a set of $k$ non-stationary regressors, $\beta $ is a $k\times 1$ vector
of the common slope parameters, and $e_{it}$ is the regression error.
Suppose $e_{it}$ is stationary and $iid$ across $i$. Then it is easy to show
that the pooled least squares estimator of $\beta $ defined by
\begin{equation}
\hat{\beta}_{LS}=\left( \sum_{i=1}^{n}\sum_{t=1}^{T}x_{it}x_{it}^{^{\prime
}}\right) ^{-1}\sum_{i=1}^{n}\sum_{t=1}^{T}x_{it}y_{it}  \label{pool-ols}
\end{equation}%
is, in general, $T$ consistent.\footnote{%
The estimator can be regarded as $\sqrt{n}T$ consistent but with a
bias of order $O(\sqrt{n}$). Up to the bias, the estimator is also
asymptotically mixed normal.} Similar to the case of time series
regression considered by Phillips and Hansen (1990), the limiting
distribution is shifted away from
zero due to an asymptotic bias induced by the long run correlation between $%
e_{it}$ and $\varepsilon _{it}$. The exception is when $x_{it}$ is strictly
exogenous, in which case the estimator is $\sqrt{n}T$ consistent. The
asymptotic bias can be estimated, and a panel fully-modified estimator can
be developed along the lines of Phillips and Hansen (1990) to achieve $\sqrt{%
n}T$ consistency and asymptotic normality.

The cross-section independence assumption is restrictive and difficult to
justify when the data under investigation are economic time series. In view
of comovements of economic variables and common shocks, we model the
cross-section dependence by imposing a factor structure on $e_{it}$. That
is,
\begin{equation*}
e_{it}=\lambda _{i}^{\prime }F_{t}+u_{it}
\end{equation*}%
where $F_{t}$ is a $r\times 1$ vector of latent common factors, $\lambda
_{i} $ is a $r\times 1$ vector of factor loadings and $u_{it}$ is the
idiosyncratic error. If  $F_t$ and $u_{it}$ are both stationary, then $e_{it}$
is also stationary. In this case, a consistent estimator of the regression
coefficients can still be obtained even when the cross-section dependence is
ignored, just like the fact that simultaneity bias is of second order in the
fixed $n$ cointegration framework. Using this property, Bai and Kao (2006)
considered a two-step fully modified estimator (2sFM). In the first step,
pooled OLS is used to obtain a consistent estimate of $\beta $. The
residuals are then used to construct a fully-modified (FM) estimator along
the line of Phillips and Hansen (1990). Essentially, nuisance parameters
induced by cross-section correlation are dealt with just like serial
correlation by suitable estimation of the long-run covariance matrices.

The 2sFM treats the I(0) common shocks as part of the error processes.
However, an alternative estimator can be developed by rewriting the
regression model as
\begin{equation}
y_{it}=x_{it}^{\prime }\beta +\lambda _{i}^{^{\prime }}F_{t}+u_{it}.
\label{f1a}
\end{equation}%
Moving $F_{t}$ from the error term to the regression function (treated as
parameters) is desirable for the following reason. 
If some components of $x_{it}$ are actually I(0), treating $F_{t}$ as part
of error process will yield an inconsistent estimate for $\beta $ when $%
F_{t} $ and $x_{it}$ are correlated. The simultaneity bias is now of the
same order as the convergence rate of the coefficient estimates on the I(0)
regressors. Estimating $\beta $ from (\ref{f1a}) with $F$ being I(0) was
suggested in Bai and Kao (2006), but its theory was not explored.

When $F_{t}$ is I(1), which is the primary focus of this paper, there is an
important difference between estimating $\beta $ from (\ref{f1a}) versus
pooled OLS in (\ref{pool-ols}) because the latter is no longer valid. More
precisely, if
\begin{equation*}
F_{t}=F_{t-1}+\eta _{t}
\end{equation*}%
then $e_{it}$ is I(1) and pooled OLS in (\ref{pool-ols}) is, in general, not
consistent. To see this, consider the following data generating process for $%
x_{it}$
\begin{equation}
x_{it}=\tau _{i}^{\prime }F_{t}+\xi _{it}  \label{xit2}
\end{equation}%
with $\xi _{it}$ being I(1) such that $\xi _{it}=\xi _{it-1}+\zeta _{it}$.
For simplicity, assume there is a single factor. It follows that $x_{it}$ is
I(1) and can be written as
(\ref{xit}) with $\varepsilon _{it}=\tau _{i}^{\prime }\eta _{t}+\zeta _{it}$%
. The pooled OLS can be written as%
\begin{equation*}
\hat{\beta}_{LS}-\beta =\frac{(\frac{1}{n}\sum_{i=1}^{n}\tau _{i}\lambda
_{i})(\frac{1}{T^{2}}\sum_{t=1}^{T}F_{t}^{2})}{\frac{1}{nT^{2}}%
\sum_{i=1}^{n}\sum_{t=1}^{T}x_{it}^{2}}+O_{p}(n^{-1/2})+O_{p}(T^{-1})
\end{equation*}%
If $\tau _{i}$ and $\lambda _{i}$ are correlated, or when they have non-zero
means, the first term on the right hand side is $O_{p}(1)$, implying
inconsistency of the pooled OLS. The best convergence rate is $\sqrt{n}$
when $x_{it}$ and $F_{t}$ are independent random walks. The problem arises
because as seen from (\ref{f1a}), we now have a panel model with
non-stationary regressors $x_{it}$ and $F_{t}$, and in which $u_{it}$ is
stationary by assumption. This means that $y_{it}$ conintegrates with $%
x_{it} $ \emph{and} $F_{t}$ with cointegrating vector $(1,\;-\beta ^{\prime
},\;\lambda _{i})$. Omitting $F_{t}$ creates a spurious regression problem.
It is worth noting that the cointegrating vector varies with $i$ because the
factor loading is unit specific. Estimation of the parameter of interest $%
\beta $ involves a new methodology because $F$ is unobservable.

In the rest of the paper, we will show how to obtain $\sqrt{n}T$ consistent
and asymptotically normal estimates of $\beta $ when the data generating
process is characterized by (\ref{f1a}) assuming that $x_{it}$ and $F_{t}$
are both I(1), and that $x_{it}$, $F_{t}$ and $u_{it}$ are potentially
correlated. We will refer to $F_{t}$ as the global stochastic trends since
they are shared by each cross-sectional unit. Hereafter, we write the
integral $\int_{0}^{1}W(s)ds\,$ as $\int W\,$ when there is no ambiguity. We
define $\Omega ^{1/2}$ to be any matrix such that $\Omega =\left( \Omega
^{1/2}\right) \left( \Omega ^{1/2}\right) ^{^{\prime }}$, and $BM\left(
\Omega \right) $ to denote Brownian motion with the covariance matrix $%
\Omega $. We use $\left\Vert A\right\Vert $ to denote $(tr(A^{\prime
}A))^{1/2}$, $\overset{d}{\longrightarrow }\,$to denote convergence in
distribution, $\overset{p}{\longrightarrow }$ to denote convergence in
probability, $\left[ x\right] $ to denote the largest integer less than or
equal to $x$. We let $M<\infty $ be a generic positive number, not depending
on $T$\ or $n$. We also define the matrix that projects onto the orthogonal
space of $z$ as $M_{z}=I_{T}-z\left( z^{\prime }z\right) ^{-1}z^{\prime }$.
We will use $\beta ^{0},$ $F_{t}^{0},$ and $\lambda _{i}^{0}$ to denote the
true common slope parameters, true common trends, and the true factor
loading coefficients. Denote $\left( n,T\right) \rightarrow \infty $ as the
joint limit. Denote $\left( n,T\right) _{seq}\rightarrow \infty $ as the
sequential limit, i.e., $T$ $\rightarrow \infty $ first and $n\rightarrow
\infty $ later. We use $MN(0,V)$ to denote a mixed normal distribution with
variance $V$.


Our analysis is based on the following assumptions.

\begin{assu}
\label{fl}Factors and Loadings:

\begin{enumerate}
\item[(a)] $E\left\Vert \lambda _{i}^{0}\right\Vert ^{4}\leq M.$ As $%
n\rightarrow \infty ,$ $\frac{1}{n}\sum_{i=1}^{n}\lambda _{i}^{0}\lambda
_{i}^{0^{\prime }}\overset{p}{\longrightarrow }\Sigma _{\lambda },$ a $%
r\times r$ diagonal matrix.

\item[(b)] $E\left\Vert \eta_t\right\Vert ^{4+\delta}\leq M$ for some $%
\delta>0$ and for all $t$. As $T\rightarrow \infty ,$ $\frac{1}{T^{2}}%
\sum_{i=1}^{n}F_{t}^{0}F_{t}^{0^{\prime }}\overset{d}{\longrightarrow }\int
B_{\eta }B_{\eta }^{^{\prime }},$ a $r\times r$ random matrix, where $%
B_{\eta }$ is a vector of Brownian motions with covariance matrix $\Omega
_{\eta }$, which is a positive definite matrix.
\end{enumerate}
\end{assu}

\begin{assu}
\label{error}Let $w_{it}=\left( u_{it},\text{ }\varepsilon _{it}^{^{\prime
}},\eta _{t}^{^{\prime }}\right) ^{^{\prime }}.$ For each $i$, $w_{it}=\Pi
_{i}(L)v_{it}=\sum_{j=0}^{\infty }\Pi _{ij}v_{it-j}$ where $v_{it}$ is
i.i.d.\ over $t $, $\sum_{j=0}^{\infty }j^{a}\left\Vert \Pi _{ij}\right\Vert
\leq M$, for some $a>1$, and $\left\vert \Pi _{i}(1)\right\vert >c> 0\,$ for
all $i$. In addition, $Ev_{it}=0$, $E(v_{it}v_{it}^{\prime })=\Sigma _{v}>0$%
, and $E\Vert v_{it}\Vert ^{8}\leq M<\infty $.
\end{assu}

\begin{assu}
\label{ser} Weak cross-sectional correlation and heterokedasticity

\begin{enumerate}
\item[(a)] $E\left( u_{it}u_{js}\right) =\sigma _{ij,ts},$ $\left\vert
\sigma _{ij,ts}\right\vert \leq \bar{\sigma }_{ij}$ for all $\left(
t,s\right) $ and $\left\vert \sigma _{ij,ts}\right\vert \leq \tau _{ts}$ for
all $\left( i,j\right) $ such that (i) $\frac{1}{n}\sum_{i,j=1}^{n}\bar{%
\sigma }_{ij}\leq M,$ (ii) $\frac{1}{T}\sum_{t,s=1}^{T}\tau _{ts}\leq M$,
and (iii) $\frac{1}{nT}\sum_{i,j,t,s=1} \left\vert \sigma
_{ij,ts}\right\vert \leq M.$

\item[(b)] For every $\left( t,s\right) ,$ $E\left\vert \frac{1}{\sqrt{n}}%
\sum_{i=1}^{n}\left[ u_{is}u_{it}-E\left( u_{is}u_{it}\right) \right]
\right\vert ^{4}\leq M. $

\item[(c)] $\frac{1}{nT^{2}}\sum_{t,s,u,v}\sum_{i,j}\left\vert cov\left(
u_{it}u_{is},u_{ju}u_{jv}\right) \right\vert \leq M $ and $\frac{1}{nT^{2}}%
\sum_{t,s}\sum_{i,j,k,l}\left\vert cov\left(
u_{it}u_{js},u_{ku}u_{ls}\right) \right\vert \leq M. $
\end{enumerate}
\end{assu}

\begin{assu}
\label{regressors}$\left\{ x_{it},F_{t}^{0}\right\} \,$are not cointegrated.
\end{assu}

Assumption \ref{fl} is standard in the panel factor literature. Assumption %
\ref{ser} allows for limited time series and cross-sectional dependence in
the error term, $u_{it}.$ Heteroskedasticity in both time series and
cross-sectional dimensions for $u_{it}$ is allowed as well. The assumption
that $\Omega _{\eta }$ is positive definite rules out cointegration among
the components of $F_{t}^{0}.$ Assumption \ref{regressors} also rules out
the cointegration between $x_{it}$ and $F_{t}^{0}.$

Assumption \ref{error} implies that a multivariate invariance principle for $%
w_{it}$ holds, i.e., the partial sum process $\frac{1}{\sqrt{T}}\sum_{t=1}^{%
\left[ T\cdot \right] }w_{it}$ satisfies:
\begin{equation*}
\frac{1}{\sqrt{T}}\sum_{t=1}^{\left[ T\cdot \right] }w_{it}\overset{d}{%
\longrightarrow }B_{i}\left( \cdot \right) =B\left( \Omega _{i}\right) \text{
as }T\rightarrow \infty \text{ for all }i,
\end{equation*}%
where
\begin{equation*}
B_{i}=\left[
\begin{array}{lll}
B_{ui} & B_{\varepsilon i}^{^{\prime }} & B_{\eta }^{^{\prime }}%
\end{array}%
\right] ^{\prime }.
\end{equation*}%
The long-run covariance matrix of $\left\{ w_{it}\right\} $ is given by
\begin{equation}
\Omega _{i}=\sum_{j=-\infty }^{\infty }E\left( w_{i0}w_{ij}^{^{\prime
}}\right) =\left[
\begin{array}{ccc}
\Omega _{ui} & \Omega _{u\varepsilon i} & \Omega _{u\eta i} \\
\Omega _{\varepsilon ui} & \Omega _{\varepsilon i} & \Omega _{\varepsilon
\eta i} \\
\Omega _{\eta ui} & \Omega _{\eta \varepsilon i} & \Omega _{\eta }%
\end{array}%
\right]  \label{omega}
\end{equation}%
are partitioned conformably with $w_{it}.$ Define the one-sided long-run
covariance
\begin{equation}
\Delta _{i}=\sum_{j=0}^{\infty }E\left( w_{i0}w_{ij}^{^{\prime }}\right) =
\left[
\begin{array}{ccc}
\Delta _{ui} & \Delta _{u\varepsilon i} & \Delta _{u\eta i} \\
\Delta _{\varepsilon ui} & \Delta _{\varepsilon i} & \Delta _{\varepsilon
\eta i} \\
\Delta _{\eta ui} & \Delta _{\eta \varepsilon i} & \Delta _{\eta }%
\end{array}%
\right] .  \label{Deltai}
\end{equation}

For future reference, it will be convenient to group elements corresponding
to $\varepsilon _{it}$ and $\eta _{t}$ taken together. Let
\begin{equation*}
B_{bi}=\left[
\begin{array}{ll}
B_{\varepsilon i}^{^{\prime }} & B_{\eta }^{^{\prime }}%
\end{array}%
\right] ^{\prime }\quad \quad \Omega _{bi}=\left[
\begin{array}{cc}
\Omega _{\varepsilon i} & \Omega _{\varepsilon \eta i} \\
\Omega _{\eta \varepsilon i} & \Omega _{\eta }%
\end{array}%
\right] .
\end{equation*}%
Then $B_{i}$ can be rewritten as
\begin{equation*}
B_{i}=\left[
\begin{array}{l}
B_{ui} \\
B_{bi}%
\end{array}%
\right] =\left[
\begin{array}{cc}
\Omega _{u.bi}^{1/2} & \Omega _{ubi}\Omega _{bi}^{-1/2} \\
0 & \Omega _{bi}^{1/2}%
\end{array}%
\right] \left[
\begin{array}{l}
V_{i} \\
W_{i}%
\end{array}%
\right]
\end{equation*}%
where $\left[
\begin{array}{ll}
V_{i} & W_{i}^{^{\prime }}%
\end{array}%
\right] ^{\prime }=BM\left( I\right) $ is a standardized Brownian motion and%
\begin{equation*}
\Omega _{u.bi}=\Omega _{ui}-\Omega _{ubi}\Omega _{bi}^{-1}\Omega _{bui}
\end{equation*}%
is the long-run conditional variance of $u_{it}$ given $(\bigtriangleup
x_{it}^{^{\prime }},\text{ }\bigtriangleup F_{t}^{^{0\prime }})^{^{\prime
}}. $ Note that $\Omega _{bi}>0$ since we assume that there is no
cointegration relationship in $(x_{it}^{^{\prime }},\text{ }F_{t}^{0^{\prime
}})^{^{\prime }}$ in Assumption \ref{regressors}.

Finally, we state an additional assumption, which is needed when deriving
the limiting distribution but is not needed for consistency of the
proposed estimators. 

\begin{assu}
\label{indepdence} The idiosyncratic errors $u_{it}$ are cross-sectionally
independent.
\end{assu}

\section{Estimation}

In this section, we first consider the problem of estimating $\beta$ when $F$
is observed. We then consider two iterative procedures that jointly estimate
$\beta$ and $F$. The procedures yield two estimators that are $\sqrt{n}T$
consistent and asymptotically normal. These estimators, denoted CupBC and
CupFM, are presented in subsections 3.2 and 3.3.

\subsection{Estimation when $F$ is observed}

The true model (\ref{f1a}) in vector form, is
\begin{equation*}
y_{i}=x_{i}\beta ^{0}+F^{0}\lambda _{i}^{0}+u_{i}
\end{equation*}%
where
\begin{equation*}
y_{i}=\left[
\begin{array}{c}
y_{i1} \\
y_{i2} \\
\vdots \\
y_{iT}%
\end{array}%
\right] ,x_{i}=\left[
\begin{array}{c}
x_{i1}^{^{\prime }} \\
x_{i2}^{^{\prime }} \\
\vdots \\
x_{iT}^{^{\prime }}%
\end{array}%
\right] ,F=\left[
\begin{array}{c}
F_{1}^{^{\prime }} \\
F_{2}^{^{\prime }} \\
\vdots \\
F_{T}^{^{\prime }}%
\end{array}%
\right] ,u_{i}=\left[
\begin{array}{c}
u_{i1} \\
u_{i2} \\
\vdots \\
u_{iT}%
\end{array}%
\right] .
\end{equation*}%
Define $\Lambda =\left( \lambda _{1},..,\lambda _{n}\right) ^{^{\prime }}$
to be an a $n\times r$ matrix. In matrix notation%
\begin{equation*}
y=X\beta ^{0}+F^{0}\Lambda ^{0^{\prime }}+u.
\end{equation*}%
Given data $y,$ $x,$ and $F^{0}$, the least squares objective function is
\begin{equation*}
S_{nT}^{0}\left( \beta ,\Lambda \right) =\sum_{i=1}^{n}\left( y-x_{i}\beta
-F^{0}\lambda _{i}\right) ^{^{\prime }}\left( y-x_{i}\beta -F^{0}\lambda
_{i}\right) .
\end{equation*}%
After concentrating out $\lambda $, the least squares estimator for $\beta $
is then%
\begin{equation*}
\tilde{\beta}_{LS}=\left( \sum_{i=1}^{n}x_{i}^{^{\prime
}}M_{F^{0}}x_{i}\right) ^{-1}\sum_{i=1}^{n}x_{i}^{^{\prime }}M_{F^{0}}y_{i}.
\end{equation*}%
The least squares estimator has the following properties.\footnote{%
The limiting distribution for $F$ being I(0) can also be obtained. Park and
Phillips (1988) provide the limiting theory with mixed I(1) and I(0)
regressors in a single equation framework.}
\begin{prop}
\label{pro1}Under Assumptions 1-5, as $\left( n,T\right) _{\func{seq}%
}\rightarrow \infty $%
\begin{equation*}
\sqrt{n}T\left( \widetilde{\beta }_{LS}-\beta ^{0}\right) -\sqrt{n}\phi
_{nT}^{0}\overset{d}{\longrightarrow }MN\left( 0,\Sigma ^{0}\right)
\end{equation*}%
where
\begin{eqnarray}
\phi _{nT}^{0} &=&\left[ \frac{1}{nT^{2}}\sum_{i=1}^{n}x_{i}^{^{\prime
}}M_{F^{0}}x_{i}\right] ^{-1}\bigg[\frac{1}{n}\sum_{i=1}^{n}\theta _{i}^{0}%
\bigg]  \label{phi} \\
\Sigma ^{0} &=&D^{-1}\left[ \underset{n\rightarrow \infty }{\lim }\frac{1}{n}%
\sum_{i=1}^{n}\Omega _{u.bi}E\left( \int Q_{i}Q_{i}^{^{\prime }}|C\right) %
\right] D^{-1},\text{ }  \label{sigma-f-known}
\end{eqnarray}%
and with $C$ being the $\sigma $-field
generated by $\left\{ F_{t}\right\},$
\begin{eqnarray*}
D &=&\underset{n\rightarrow \infty }{\lim }\frac{1}{n}\sum_{i=1}^{n}E\left(
\int Q_{i}Q_{i}^{\prime }|C\right) \\
Q_{i} &=&B_{\varepsilon i}-\left( \int B_{\varepsilon i}B_{\eta }^{^{\prime
}}\right) \left( \int B_{\eta }B_{\eta }^{^{\prime }}\right) ^{-1}B_{\eta },
\\
\theta _{i}^{0} &=&\frac{1}{T}x_{i}^{^{\prime }}M_{F^{0}}\Delta b_{i}\Omega
_{bi}^{-1}\Omega _{bui}+\left( \Delta _{\varepsilon ui}^{+}-\delta
_{i}^{^{0\prime }}\Delta _{\eta u}^{+}\right) , \\
\delta _{i}^{0} &=&(F^{0^{\prime }}F^{0})^{-1}F^{0^{\prime }}x_{i},\quad
\Delta b_{i}=(%
\begin{array}{cc}
\Delta x_{i} & \Delta F^{0}%
\end{array}%
)
\end{eqnarray*}%
\begin{equation*}
\Delta _{bui}^{+}=\left(
\begin{array}{c}
\Delta _{\varepsilon ui}^{+} \\
\Delta _{\eta u}^{+}%
\end{array}%
\right) =\left(
\begin{array}{ll}
\Delta _{bui} & \Delta _{bi}%
\end{array}%
\right) \left(
\begin{array}{c}
I_{k} \\
-\Omega _{bi}^{-1}\Omega _{bui}%
\end{array}%
\right) =\Delta _{bui}-\Delta _{bi}\Omega _{bi}^{-1}\Omega _{bui}.
\end{equation*}
\end{prop}

The estimator is $\sqrt{n}T$ consistent if $\phi^0 _{nT}=0$, which occurs
when $x_{it}$ is strictly exogenous. Otherwise, the estimator is $T$
consistent as there is an asymptotic bias given by the term $\sqrt{n}
\phi^0_{nT}$. This is an average of individual biases that are data specific
as seen from the definition of $\theta_i^0$. The individual biases arise
from the contemporaneous and low frequency correlations between the
regression error and the innovations of the I(1) regressors as given by
terms such as $\Omega _{bui} $ and $\Delta_{bui}.$

To estimate the bias, we need to consistently estimate the nuisance
parameters. We use a kernel estimator. Let
\begin{eqnarray*}
\widehat{\Omega }_{i} &=&\sum_{j=T+1}^{T-1}\omega\left( \frac{j}{K}\right)
\widehat{\Gamma }_{i}\left( j\right) , \\
\widehat{\Delta }_{i} &=&\sum_{j=0}^{T-1}\omega\left( \frac{j}{K}\right)
\widehat{\Gamma }_{i}\left( j\right) \\
\widehat{\Gamma }_{i}\left( j\right) &=&\frac{1}{T}\sum_{t=1}^{T-j}\widehat{w%
}_{it+j}\widehat{w}_{it}^{^{\prime }}.
\end{eqnarray*}%
where $\hat w_{it}=(\hat u_{it}, \Delta x_{it}^{\prime }, \Delta
F_t^{0\prime} )^{\prime }$. To state the asymptotic theory for the
bias-corrected estimator, we need the following assumption, as used in Moon
and Perron (2004):

\begin{assu}
\label{band}

\begin{enumerate}
\item[(a)] $\underset{ n,T \rightarrow \infty } \liminf \, \, ( \log T/\log
n) >1.$

\item[(b)] the kernel function $\omega\left( \cdot \right) :\mathrm{R}%
\rightarrow \left[ -1,1\right] $ satisfies (i) $\omega\left( 0\right) =1,$ $%
\omega\left( x\right) =\omega\left( -x\right) ,$ (ii) $\int_{-1}^{1}\omega%
\left( x\right) ^{2}dx<\infty $ and with Parzen's exponent $q\in \left(
0,\infty \right) $ such that $\lim \frac{1-\omega\left( x\right) }{%
\left\vert x\right\vert ^{q}}<\infty . $

\item[(c)] The bandwidth parameter $K$ satisfies $K\backsim n^{b}$ and $%
\frac{1}{2q}<b<\lim \inf \frac{\log T}{\log n}-1.$
\end{enumerate}
\end{assu}

Let
\begin{equation*}
\hat{\phi}_{nT}^{0}=\left[ \frac{1}{nT^{2}}\sum_{i=1}^{n}x_{i}^{^{\prime
}}M_{F^{0}}x_{i}\right] ^{-1}\hat{\theta}^{n}
\end{equation*}%
where $\hat{\theta}^{n}=\frac{1}{n}\sum_{i=1}^{n}\hat{\theta}_{i}$, $\hat{%
\theta}_{i}$ is a consistent estimate of $\theta _{i}^{0}$. The resulting
bias-corrected estimator is
\begin{equation}
\tilde{\beta}_{LSBC}=\tilde{\beta}_{LS}-\frac{1}{T}\hat{\phi}_{nT}^{0}.
\label{lsbc}
\end{equation}%
This estimator can alternatively be written as
\begin{equation}
\widetilde{\beta }_{LSFM}=\left( \sum_{i=1}^{n}x_{i}^{^{\prime
}}M_{F^{0}}x_{i}\right) ^{-1}\sum_{i=1}^{n}\left( x_{i}^{^{\prime }}M_{F^{0}}%
\tilde{y}_{i}^{+}-T\left( \tilde{\Delta}_{\varepsilon ui}^{+}-\delta
_{i}^{0\prime }\tilde{\Delta}_{\eta u}^{+}\right) \right)  \label{lsfm}
\end{equation}%
where $\tilde{y}^{+}$ and $\tilde{\Delta}^{+}$ are consistent estimates of $%
y^{+}$ and $\Delta ^{+}$ etc, with
\begin{equation*}
y_{it}^{+}=y_{it}-\Omega _{ubi}\Omega _{bi}^{-1}\left(
\begin{array}{c}
\Delta x_{it} \\
\Delta F_{t}^{0}%
\end{array}%
\right) \quad \quad u_{it}^{+}=u_{it}-\Omega _{ubi}\Omega _{bi}^{-1}\left(
\begin{array}{c}
\Delta x_{it} \\
\Delta F_{t}^{0}%
\end{array}%
\right)
\end{equation*}%
Viewed in this light, the bias-corrected estimator is also a panel
fully-modified estimator in the spirit of Phillips and Hansen (1990), and is
the reason why the estimator is also labeled $\hat\beta_{LSFM}$. It is not
difficult to verify that $\hat{\beta}_{LSBC}$ and $\hat{\beta}_{LSFM}$ are
identical. Panel fully modified estimators were also considered by Phillips
and Moon (1999) and Bai and Kao (2006). Here, we extend those analysis to
allow for common stochastic trends. By construction $u_{it}^{+}$ has a zero
long-run covariance with $(%
\begin{array}{cc}
\Delta x_{it}^{\prime } & \Delta F_{t}^{0\prime }%
\end{array}%
)^{^{\prime }}$ and hence the endogeneity can be removed. Furthermore,
nuisance parameters arising from the low frequency correlation of the errors
are summarized in $\Delta _{bui}^{+}$.

\begin{prop}
\label{pro2}Let $\widetilde{\beta }_{LSFM}$ be defined by (\ref{lsfm}).
Under Assumptions 1-6, as $\left( n,T\right) _{seq}\rightarrow \infty $

\begin{equation*}
\sqrt{n}T(\widetilde{\beta }_{LSFM}-\beta ^{0})\overset{d}{\longrightarrow }%
MN\left( 0,\Sigma ^{0}\right) .
\end{equation*}
\end{prop}

In small scale cointegrated systems, cointegrated vectors are $T$
consistent, and this fast rate of convergence is already accelerated
relative to the case of stationary regressions, which is $\sqrt{T}$. Here in
a panel data context with observed global stochastic trends, the estimates
converge to the true values at an even faster rate of $\sqrt{n}T$ and the
limiting distributions are normal. To take advantage of this fast
convergence rate made possible by large panels, we need to deal with the
fact that $F^0 $ is not observed. This problem is considered in the next two
subsections.

\subsection{Unobserved $F^0$ and the Cup Estimator}

The LSFM considered above is a linear estimator and can be obtained if $F^0$
is observed. When $F^0$ is not observed, the previous estimator is infeasible.
Recall that least squares estimator that ignores $F$ is, in general,
inconsistent. In this section, we consider estimating $F$ along with $\beta$
and $\Lambda$ by minimizing the objective function
\begin{equation}
S_{nT}\left( \beta ,F,\Lambda \right) =\sum_{i=1}^{n}\left( y-x_{i}\beta
-F\lambda _{i}\right) ^{^{\prime }}\left( y-x_{i}\beta -F\lambda _{i}\right)
\label{ob}
\end{equation}%
subject to the constraint $T^{-2}F^{^{\prime }}F =I_{r}$ and $\Lambda
^{^{\prime }}\Lambda $ is positive definite. The least squares estimator for $%
\beta $ for a given $F$ is
\begin{equation*}
\widehat{\beta }=\left( \sum_{i=1}^{n}x_{i}^{^{\prime }}M_{F}x_{i}\right)
^{-1}\sum_{i=1}^{n}x_{i}^{^{\prime }}M_{F}y_{i}.
\end{equation*}%
Define
\begin{eqnarray*}
w_{i}&=&y_{i}-x_{i}\beta \\
&=&F\lambda _{i}+u_{i}.
\end{eqnarray*}
Notice that given $\beta ,$ $w_i$ has a pure factor structure. Let $W=\left(
w_{i},...,w_{n}\right) $ be a $T\times n$ matrix. We can rewrite the
objective function (\ref{ob}) as $tr[( W-F\Lambda ^{^{\prime }}) (
W-F\Lambda ^{\prime}) ^{\prime}] . $ If we concentrate out $\Lambda
=W^{^{\prime }}F\left( F^{^{\prime }}F\right) ^{-1}=T^{-2} W^{^{\prime }}F$,
we have the concentrated objective function:
\begin{equation}
tr\left( W^{^{\prime }}M_{F}W\right) =tr\left( W^{^{\prime }}W\right)
-tr\left( F^{^{\prime }}WW^{^{\prime }}F /T^2 \right) .  \label{ob2}
\end{equation}
Since the first term does not depend on $F$, minimizing (\ref{ob2}) with
respect to $F$ is equivalent to maximizing $tr\left( T^{-2}F^{^{\prime
}}WW^{^{\prime }}F\right) $ subject to the constraint $T^{-2}F^{\prime
}F=I_{r}.$ The solution, denoted $\widehat{F}$, is a matrix of the first $r$
eigenvectors (multiplied by $T$) of the matrix $\frac{1}{nT^{2}}%
\sum_{i=1}^{n}\left( y_{i}-x_{i}\beta \right) \left( y_{i}-x_{i}\beta
\right) ^{^{\prime }}$.

Although $F$ is not observed when estimating $\beta $, and similarly, $\beta
$ is not observed when estimating $F$, we can replace the unobserved
quantities by initial estimates and iterate until convergence. Such a
solution is more easily seen if we rewrite the left hand side of (\ref{ob2})
with $y-x\beta $ substituting in for $W$. Define
\begin{equation*}
S_{nT}\left( \beta ,F\right) =\frac{1}{nT^{2}}\sum_{i=1}^{n}\left(
y_{i}-x_{i}\beta \right) ^{^{\prime }}M_{F}\left( y_{i}-x_{i}\beta
\right).
\end{equation*}%
The continuous updated estimator (Cup) for $\left( \beta ,F\right) $ is
defined as
\begin{equation*}
\left( \widehat{\beta }_{Cup},\widehat{F}_{Cup}\right) =\underset{\beta ,F}{%
\text{ argmin}} \, S_{nT}\left( \beta ,F\right) .
\end{equation*}%
More precisely, $(\widehat{\beta }_{Cup},\widehat{F}_{Cup})$ is the solution
to the following two nonlinear equations%
\begin{eqnarray}
\widehat{\beta } &=&\left( \sum_{i=1}^{n}x_{i}^{^{\prime }}M_{\widehat{F}%
}x_{i}\right) ^{-1}\sum_{i=1}^{n}x_{i}^{^{\prime }}M_{\widehat{F}}y_{i}
\label{ls} \\
\widehat{F}V_{nT} &=&\left[ \frac{1}{nT^{2}}\sum_{i=1}^{n}\left( y_{i}-x_{i}%
\widehat{\beta }\right) \left( y_{i}-x_{i}\widehat{\beta }\right) ^{^{\prime
}}\right] \widehat{F}  \label{ls1}
\end{eqnarray}%
where $M_{\widehat{F}}=I_{T}-T^{-2}\widehat{F}\widehat{F}^{\prime }$ since $%
\widehat{F}^{\prime }\widehat{F}/T^{2}=I_{r}$, and $V_{nT}$ is a diagonal
matrix consisting of the $r$ largest eigenvalues of the matrix inside the
brackets, arranged in decreasing order. Note that the estimator is obtained
by iteratively solving for $\hat{\beta}$ and $\hat{F}$ using (\ref{ls}) and (%
\ref{ls1}). It is a non-linear estimator even though linear least squares
estimation is involved at each iteration. An estimate of $\Lambda $ can be
obtained as:
\begin{equation*}
\widehat{\Lambda }=T^{-2}\widehat{F}^{^{\prime }}\left( Y-X\widehat{\beta }%
\right) .
\end{equation*}%
The triplet $\left( \widehat{\beta },\widehat{F},\widehat{\Lambda }\right) $
jointly minimizes the objective function (\ref{ob}).

The estimator $\hat \beta_{Cup}$ is consistent for $\beta$. We state this
result in the following proposition.

\begin{prop}
\label{consistency} Under Assumptions \ref{fl}-\ref{regressors} and as $\left(
n,T\right) \rightarrow \infty $,
\begin{equation*}
\widehat{\beta }_{Cup}\overset{p}{\longrightarrow }\beta ^{0}.
\end{equation*}
\end{prop}

We now turn to the asymptotic representation of $\hat\beta_{Cup}$.

\begin{prop}
\label{known-f}Suppose Assumptions \ref{fl}-\ref{regressors} hold and $%
\left( n,T\right) \rightarrow \infty $. Then
\begin{equation*}
\sqrt{n}T\left( \widehat{\beta }_{Cup}-\beta ^{0}\right) =D\left(
F^{0}\right) ^{-1}\left[ \frac{1}{\sqrt{n}T}\sum_{i=1}^{n}\left(
x_{i}^{^{\prime }}M_{F^{0}}-\frac{1}{n}\sum_{k=1}^{n}a_{ik}x_{i}^{^{\prime
}}M_{F^{0}}\right) u_{i}\right] +o_{p}\left( 1\right) ,
\end{equation*}%
where $a_{ik}=\lambda _{i}^{\prime }\left( \frac{\Lambda ^{^{\prime
}}\Lambda }{n}\right) ^{-1}\lambda _{k}$, $D\left( F^{0}\right) =\frac{1}{%
nT^{2}}\sum_{i=1}^{n}Z_{i}^{\prime }Z_{i}$ and $Z_{i}=M_{F^{0}}x_{i}-\frac{1%
}{n}\sum_{k=1}^{n}M_{F^{0}}x_{k}a_{ik}.$
\end{prop}

In comparison with the pooled least squares estimator for the case of known $%
F^{0}$, estimation of the stochastic trends clearly affects the limiting
behavior of the estimator. The term involving $a_{ik}$ is due to the
estimation of $F$. This effect is carried over to the limiting distribution
and to the asymptotic bias, as we now proceed to show. Let $\bar{w}_{it}=(
u_{it},\text{ }\Delta \bar{x}_{i}^{\prime },\eta _{t}^{\prime })^{\prime }$
where $\bar{x}_{i}=x_{i}-\frac{1}{n}\sum_{k=1}^{n}x_{k}a_{ik}$. For the rest
of the paper, we use bar to denote those long run covariance matrices
(including one sided and conditional covariances and so on) generated from $%
\bar{w}_{it}$ instead of $w_{it}.$ Thus, $\bar{\Omega}_{i}$ is the long run
covariance matrix of $\bar{w}_{it}$ as in (\ref{omega}), and define $\bar{%
\Delta}_{i}$ is the one-sided covariance matrix of $\bar w_{it}$. These
quantities depend on $n$, but this dependence is suppressed for notional
simplicity.

Because the right hand side of the representation does not depend on
estimated quantities, it is not difficult to derive the limiting
distribution of $\hat \beta_{Cup}$, even allowing for cross-sectional
correlation in $u_{it}$. However, estimating the resulting nuisance
parameters would be more difficult. Thus, although consistency of the Cup
estimator does not require the cross-section independence of $u_{it}$, our
asymptotic distribution for $\hat \beta_{Cup}$ is derived with Assumption %
\ref{indepdence} imposed.

\begin{thm}
\label{theo1}Suppose that Assumptions \ref{fl}-5 hold. Let $\hat{\beta}%
_{Cup} $ be obtained by iteratively updating (\ref{ls}) and (\ref{ls1}). As $%
\left( n,T\right) _{seq}\rightarrow \infty $, we have
\begin{equation*}
\sqrt{n}T\left( \widehat{\beta }_{Cup}-\beta \right) -\sqrt{n}\phi _{nT}%
\overset{d}{\longrightarrow }MN\left( 0,\Sigma \right)
\end{equation*}%
where
\begin{equation*}
\phi _{nT}=\left[ \frac{1}{nT^{2}}\sum_{i=1}^{n}Z_{i}^{^{\prime }}Z_{i}%
\right] ^{-1}\Big(\frac{1}{n}\sum_{i=1}^{n}\theta _{i}\Big)
\end{equation*}%
\begin{equation*}
\theta _{i}=\frac{1}{T}Z_{i}^{^{\prime }}\Delta \bar{b}_{i}\bar{\Omega}%
_{bi}^{-1}\bar{\Omega}_{bui}+\left( \bar{\Delta}_{\varepsilon ui}^{+}-\bar{%
\delta}_{i}^{^{\prime }}\bar{\Delta}_{\eta u}^{+}\right) ,
\end{equation*}%
\begin{equation}
\Sigma =D_{Z}^{-1}\left[ \underset{n\rightarrow \infty }{\lim }\frac{1}{n}%
\sum_{i=1}^{n}\bar{\Omega}_{u.bi}E\left( \int R_{ni}R_{ni}^{^{\prime
}}|C\right) \right] D_{Z}^{-1},\text{ }D_{Z}=\text{ }\underset{n\rightarrow
\infty }{\lim }\frac{1}{n}\sum_{i=1}^{n}E\left( \int R_{ni}R_{ni}^{^{\prime
}}|C\right) ,  \label{sigma-f-unknown}
\end{equation}%
%
%
%
%
%
%
%
%
%
%
%
%
%
%
%
%
%
%
%
%
%
\begin{equation*}
R_{ni}=Q_{i}-\frac{1}{n}\sum_{k=1}^{n}Q_{k}a_{ik},
\end{equation*}%
\begin{eqnarray*}
\Delta \bar{b}_{i} &=&\left(
\begin{array}{cc}
\Delta \bar{x}_{i} & \Delta F^{0}%
\end{array}%
\right) , \\
\bar{x}_{i} &=&x_{i}-\frac{1}{n}\sum_{k=1}^{n}x_{k}a_{ik}, \\
\bar{\delta}_{i} &=&\delta _{i}-\frac{1}{n}\sum_{k=1}^{n}\delta _{k}a_{ik}.
\end{eqnarray*}
\end{thm}


Theorem 1 establishes the large sample properties of the Cup
estimator. The Cup estimator is $\sqrt{n}T$ consistent provided that
$\phi_{nT}=0$, which occurs when $x_{it}$ and $F_t$ are exogenous.
Since $\phi _{nT}=O_{p}(1)$, the Cup estimator is at least $T$
consistent. This is in contrast with pooled OLS in Section 2, where
it was shown to be inconsistent in general. Nevertheless, as in the
case when $F$ is observed, the Cup estimator has an asymptotic bias
and thus the limiting distribution is not centered around zero.
There is an extra bias term (the term involving $a_{ik})$ that
arises from having to estimate $F_t$.  In consequence, the bias is
now a function of terms not present in Proposition 1, which is valid
when $F_{t}$ is observed. We now consider removing the bias by
constructing a
consistent estimate of $\phi _{nT}$. This can be obtained upon replacing $%
F^{0}$, $\Delta \bar b_{i},$ $\bar{\Omega }_{bi},$ $\bar{\Omega }_{bui},$ $%
\bar{\Delta }_{\varepsilon ui}^+,$ $\bar{\Delta }_{\eta u}^+$ by their
consistent estimates.

We consider two fully-modified estimators. The first one directly corrects
the bias of $\hat \beta_{Cup}$, and is denoted by $\hat \beta_{CupBC}$. The
second one will be considered in the next subsection, where correction is
made during each iteration, and will be denoted by $\hat \beta_{CupFM}$.

Consider
\begin{eqnarray*}
\widehat{\bar{\Omega}}_{i} &=&\sum_{j=T+1}^{T-1}\omega\left( \frac{j}{K}%
\right) \widehat{\Gamma }_{i}\left( j\right) , \\
\widehat{\bar{\Delta}}_{i} &=&\sum_{j=0}^{T-1}\omega\left( \frac{j}{K}%
\right) \widehat{\Gamma }_{i}\left( j\right) \\
\widehat{\Gamma }_{i}\left( j\right) &=&\frac{1}{T}\sum_{t=1}^{T-j}\widehat{%
\bar{w}}_{it+j}\widehat{\bar{w}}_{it}^{^{\prime }}.
\end{eqnarray*}%
where
\begin{equation*}
\widehat{\bar{w}}_{it}=(\hat{u}_{it},\Delta \hat{\bar{x}}_{it}^{\prime
},\Delta \hat{F}_{t}^{\prime })^{\prime }\quad \mathrm{with~}\Delta {\hat{%
\bar{x}}_{it}}=\Delta x_{it}-\frac{1}{n}\sum_{k=1}^{n}\Delta x_{kt}\hat{a}%
_{ik}
\end{equation*}%
The bias-corrected Cup estimator is defined as
\begin{equation*}
\hat{\beta}_{CupBC}=\hat{\beta}_{Cup}-\frac{1}{T}\hat{\phi}_{nT}
\end{equation*}%
where
\begin{eqnarray*}
\widehat{\phi }_{nT} &=&\left[ \frac{1}{nT^{2}}\sum_{i=1}^{n}\widehat{Z}%
_{i}^{^{\prime }}\widehat{Z}_{i}\right] ^{-1} \Big(\frac{1}{n}\sum_{i=1}^{n}%
\hat{\theta}_{i}\Big) \\
\hat \theta_i &=&\widehat{Z}_{i}^{^{\prime }}\Delta \widehat{\bar{b}}_{i}%
\hat{\bar{\Omega}}_{bi}^{-1}\hat{\bar{\Omega}}_{bui}+\left(
\begin{array}{cc}
\widehat{\bar{\Delta}}_{\varepsilon ui}^{+} & -\hat{\bar{\delta}}%
_{i}^{^{\prime }}%
\end{array}%
\widehat{\Delta }_{\eta u}^{+}\right) , \\
\hat{\bar \delta}_{i} &=&\left( \widehat{F}^{\prime }\widehat{F}\right) ^{-1}%
\widehat{F}^{\prime }\widehat{\bar{x}}_{i}\quad \quad \Delta \widehat{b}%
_{i}=\left(
\begin{array}{cc}
\Delta \hat {\bar{x}}_{i} & \Delta \widehat{F}%
\end{array}%
\right) \\
\widehat{\bar{x}}_{i} &=&x_{i}-\frac{1}{n}\sum_{k=1}^{n}x_{kt}\widehat{a}%
_{ik},\quad \quad \widehat{a}_{ik}=\widehat{\lambda }_{i}^{^{\prime }}\left(
\widehat{\Lambda }^{^{\prime }}\widehat{\Lambda }/n\right) ^{-1}\widehat{%
\lambda }_{k}.
\end{eqnarray*}

\begin{thm}
Suppose Assumptions \ref{fl}-\ref{band} hold. Then as $\left( n,T\right)
_{seq}\rightarrow \infty$,
\begin{equation*}
\sqrt{n}T\left( \widehat{\beta }_{CupBC}-\beta ^{0}\right) \overset{d}{%
\longrightarrow }MN\left( 0,\Sigma \right) .
\end{equation*}
\end{thm}

The CupBC is $\sqrt{n}T$ consistent with a limiting distribution
that is centered at zero. This type of bias correction approach is
also used in Hahn and Kuersteiner (2002), for example, and is not
uncommon in panel data analysis. Because the bias-corrected
estimator is $\sqrt{n}T $ and has a normal limit distribution, the
usual $t$ and Wald tests can be used for inference. Note that the
limiting distribution is different from that of the infeasible LSBC
estimator, which coincides with LSFM and whose asymptotic variance
is $\Sigma^0$ instead of $\Sigma$. Thus, the estimation of $F$
affects the asymptotic distribution of the estimator. As in the case
when $F$ is observed, the bias corrected estimator can be rewritten
as a fully modified estimator. Such a fully-modified estimator is
now discussed.

\subsection{ A Fully Modified Cup Estimator}

The CupBC just considered is constructed by estimating the asymptotic bias
of $\hat \beta_{Cup}$, and then subtracting it from $\hat \beta_{Cup}$. In
this subsection, we consider a different fully-modified estimator, denoted
by $\hat \beta_{CupFM}$. Let
\begin{eqnarray*}
y_{it}^{+} &=&y_{it}-\widehat{\bar \Omega }_{ubi}\widehat{\bar \Omega }%
_{bi}^{-1}\left(
\begin{array}{c}
\Delta \widehat{\bar{x}}_{it} \\
\Delta \widehat{F}_{t}%
\end{array}%
\right) \\
\hat{\bar \delta}_{i} &=&\left( \widehat{F}^{\prime }\widehat{F}\right) ^{-1}%
\widehat{F}^{\prime }\widehat{\bar{x}}_{i}
\end{eqnarray*}%
where $\widehat{\bar \Omega }_{ubi},$ $\widehat{\bar \Omega }_{bi},$ and $%
\widehat{\bar \Delta }_{bui}$ are estimates of $\bar{\Omega }_{ubi},\bar{%
\Omega }_{bi}$ and $\bar{\Delta }_{bui}$, respectively. Recall that $%
\widehat{\beta }_{Cup}$ is obtained by jointly solving (\ref{ls}) and (\ref%
{ls1}). Consider replacing these equations by the following:
\begin{eqnarray}
\widehat{\beta }_{CupFM} &=&\left( \sum_{i=1}^{n}x_{i}^{^{\prime }}M_{%
\widehat{F}}x_{i}\right) ^{-1}\sum_{i=1}^{n}\left( x_{i}^{^{\prime }}M_{%
\widehat{F}}y_{i}^{+}-T\left( \widehat{\bar \Delta }_{\varepsilon ui}^{+}-%
\widehat{\bar \delta }_{i}^{^{\prime }}\widehat{\bar \Delta }_{\eta
u}^{+}\right) \right)  \label{cupfm} \\
\widehat{F}V_{nT} &=&\left[ \frac{1}{nT^{2}}\sum_{i=1}^{n}\left( y_{i}-x_{i}%
\widehat{\beta }_{CupFM}\right) \left( y_{i}-x_{i}\widehat{\beta }%
_{CupFM}\right) ^{^{\prime }}\right] \widehat{F}.  \label{cupfm1}
\end{eqnarray}%
Like the FM estimator of Phillips and Hansen (1990), the corrections are
made to the data to remove serial correlation and endogeneity. The CupFM
estimator for $(\beta ,F)$ is obtained by iteratively solving (\ref{cupfm})
and (\ref{cupfm1}). Thus correction to endogeneity and serial correlation is
made during each iteration.

\begin{thm}
\label{fesi-fm}Suppose Assumptions \ref{fl}-6 hold. Then as $\left(
n,T\right) _{seq}\rightarrow \infty ,$
\begin{equation*}
\sqrt{n}T\left( \widehat{\beta }_{CupFM}-\beta ^{0}\right) \overset{d}{%
\longrightarrow }MN\left( 0,\Sigma \right) ,
\end{equation*}%
where $\Sigma $ is given in (\ref{sigma-f-unknown})
\end{thm}

The CupFM and CupBC have the same asymptotic distribution, but they are
constructed differently. The estimator $\widehat{\beta }_{CupBC}$ does the
bias correction only once, i.e., at the final stage of the iteration, and $%
\widehat{\beta }_{CupFM}$ does the correction at every iteration.
The situation is different from the case of known $F$, in which the
bias-corrected estimator and the fully-modified estimator are
identical due to the absence of iteration. Again, because of the
mixture of normality, hypothesis testing on $\beta$ can proceed with
the usual t or chi square distributions.

Kapetanios, Pesaran, and Yamagata (2006) suggest an alternative
estimation procedure based  on Pesaran (2006). The  model is
augmented with additional regressors $\bar y_t$ and $\bar x_t$,
which are cross-sectional averages of  $y_{it}$ and $x_{it}$. These
averages are used as proxy for $F_t$. The estimator for the slope
parameter $\beta$ is shown to be  $\sqrt{n}$ consistent,  but fully
modified estimator is not considered.

While the focus is on estimating the slope parameters $\beta$, the
global stochastic trends $F$ are also of interest.  We state this
result as a proposition:
\begin{prop}
\label{Fhat} Let $\hat F$ be the solution of (\ref{cupfm1}). Under
Assumptions of 1-4, we have
\begin{equation*}
\frac 1 T \sum_{t=1}^T \| \hat F_t -H F_t^0 \|^2 =O_p(\frac 1 n) +
O_p( \frac 1 {T^2})
\end{equation*}
where $H$ is an $r\times r$ invertible matrix.
\end{prop}

Thus, we can estimate the true global stochastic trends up to a rotation.
This is the same rate as in Bai (2004, Lemma B.1), where the regressor $%
x_{it}$ is absent. Similarly, the factor loadings $\lambda_i$ are estimated
with the same rate of convergence as in Bai (2004).

Thus far, our analysis  assumes that the number of stochastic
trends, $r$, is known. If this is not the case, $r$ can be
consistently estimated using the information criterion function
developed in Bai and Ng (2002). In particular, let
\begin{equation*}
\widehat{r}=\arg \underset{1\leq r\leq r_{\max }}{\min } IC\left(
r\right)
\end{equation*}%
where $r\leq r_{\max }$, $r_{\max }$ is a bounded integer and
\begin{equation*}
IC(r)=\log \widehat{\sigma }^{2}\left( r\right) +r g_{nT}
\end{equation*}%
where $g_{nT}\rightarrow 0$ as $n,T\rightarrow \infty $ and $\min
[n,T]g_{nT}\rightarrow \infty $. For example, $g_{nT}$ can be $\log
(a_{nT})/a_{nT}$, with $a_{nT}=\frac{nT}{n+T}$. Then $P(\hat
r=r)\rightarrow 1$ as $n,T\rightarrow\infty$. This criterion
estimates the total number of factors, including I(0) factors. To
estimate the number of I(1) factors only, the criterion in Bai
(2004) can be used. Ignoring the I(0) factors still lead to
consistent estimation of $\beta$. However, our
distribution theory assumes cross-sectional independence for the
idiosyncratic errors;  lumping I(0) factors with the regression
errors will violate this assumption. This suggests the use of
 Bai-Ng criterion.




\section{Further issues}

The preceding analysis assumes that there are no deterministic
components and that the regressors and the common factors are all
I(1) without drifts. This section considers construction of the
estimator when these restrictions are relaxed. It will be shown that
when there are deterministic components, we can apply the same
estimation procedure to the demeaned or detrended series, and the
Brownian motion processes in the limiting distribution are replaced
by the demeaned and/or detrended versions. Furthermore, the
procedure is robust to the presence of mixed I(1)/I(0) regressors
and/or factors. Of course, the convergence rates for I(0) and I(1)
regressors will be different, but asymptotic mixed normality and the
construction of test statistics (and their limiting distribution) do
not depend on the convergence rate. Finally, we also discuss the
issue of testing cross-sectional independence.

\subsection{Incidental trends}

The Cup estimator can be easily extended to models with incidental trends,
\begin{equation}  \label{trend}
y_{it}= \alpha_i + \rho_i t + x_{it}^{\prime }\beta + \lambda_i^{\prime }F_t
+ u_{it}.
\end{equation}
In the intercept only case ($\rho_i=0$, for all $i$), we define the
projection matrix
\begin{equation*}
M_T= I_T -\iota_T \iota_T^{\prime }/T
\end{equation*}
where $\iota_T$ is a vector of 1's. When a linear trend is also included in
the estimation, we define $M_T$ to be the projection matrix orthogonal to $%
\iota_T$ and to the linear trend. Then
\begin{equation*}
M_T y_i = M_T x_i \beta + M_T F \lambda_i +M_T u_i,
\end{equation*}
or
\begin{equation*}
\dot y_i =\dot x_i \beta + \dot F_t \lambda_i + \dot u_i
\end{equation*}
where the dotted variables are demeaned and/or detrended versions. The
estimation procedure for the cup estimator is identical to that of Section
3, except that we use dotted variables.

With the intercept only case, the construction of FM estimator is also the
same as before. Theorems 1-3 hold with the following modification for the
limiting distribution. The random processes $B_{\varepsilon,i}$ and $B_\eta$
in $Q_i$ are replaced by the demeaned Brownian motions.

When linear trends are allowed, $\Delta x_{it}$ is now replaced by
$\hat \varepsilon_{it}  =\Delta x_{it}-%
\overline{\Delta x_i}$, which is detrended residual of
$x_{it}$. But since $\dot x_i$ is already a detrended series, and $\hat F$
is also asymptotically detrended (since it is estimating $\dot F$), $\Delta
\dot x_{it}$ and $\Delta \hat F_t$ are also estimating the detrended
residuals. Thus we can simply apply the same procedure prescribed in Section
3 with the dotted variables. The limiting distribution in Theorem 2 and
consequently in Theorem 3 is modified upon replacing the random processes $%
B_{\varepsilon i}$ and $B_\eta$  by the demeaned and detrended
Brownian motions.\footnote{Alternatively, we can use
$\hat \varepsilon_{it} -\frac 1 n \sum_{k=1}^n \hat \varepsilon_{kt} \hat
a_{ik}$ in place of $\Delta \hat {\bar x}_{it}$ in Section 3. Similarly,
 we use $\hat \eta_t=
\Delta \hat F_t -\overline{\Delta \hat F}$ in place of $\Delta \hat F_t$.}
The test statistics ($t$ and $\chi ^{2}$) have standard
asymptotic distribution, not depending on whether the underlying Brownian
motion is demeaned or detrended.

When linear trends are included in the estimation, the limiting distribution
is invariant to whether or not $y_{it}$, $x_{it}$ and $F_t$ contain a linear
trend. Now suppose that these variables do contain a linear trend (drifted
random walks). With deterministic cointegration holding (i.e., cointegrating
vector eliminates the trends), the estimated $\beta$ will have a faster
convergence rate when a separate linear trend is not included in the
estimation. But we do not consider this case. Interested readers are
referred to Hansen (1992).

\subsection{Mixed I(0)/I(1) Regressors and Common Shocks}

So far, we have considered estimation of panel cointegration models when all
the regressors and common shocks are I(1). There are no stationary
regressors or stationary common shocks. In this section we suggest that the
results are robust to mixed I(1)/I(0) regressors and mixed I(1)/I(0) common
shocks. Below, we sketch the arguments for the LS estimator assuming the
factors are observed. If they are not observed, the limiting distribution is
different, but the idea of argument is the same.

Recall that the LS estimator is $\hat\beta_{LS}=(\sum_{i=1}^n x_i^\prime
M_{F^0} x_i)^{-1} \sum_{i=1}^n x_i^\prime M_{F^0} y_i$. The term
\begin{equation*}
M_{F^{0}}x_{i}=( I_{T}-F^{0}( F^{0^{\prime }}F^{0}) ^{-1}F^{0^{\prime }})
x_{i}=x_{i}-F^{0}\delta _{i}
\end{equation*}
with $\delta _{i}=( F^{0^{\prime }}F^{0}) ^{-1}F^{0^{\prime }}x_{i}$ plays
an important role in the properties of the LS. When $x_{it}$ and $F_t$ are
I(1), $\delta_i=O_p(1)$ and thus
\begin{equation*}
\frac{(M_{F^0}x_i)_t}{\sqrt{T}} =\frac{x_{it}}{\sqrt{T}}-\frac{%
\delta_i^{\prime }F_t^0}{\sqrt{T}}=O_p(1).
\end{equation*}
We now consider this term under mixed I(1) and I(0) assumptions.

\paragraph{I(1) Regressors, I(0) Factors.}

Suppose all regressors are I(1) and all common shocks are I(0). With I(0)
factors, we have $T^{-1}F^{0\prime }F^{0}\overset{p}{\longrightarrow }\Sigma
_{F}=O_{p}\left( 1\right) $. Thus
\begin{equation*}
\delta _{i}=\left( T^{-1}F^{0^{\prime }}F^{0}\right) ^{-1}\frac{1}{T}%
\sum_{t=1}^{T}F_{t}^{0}x_{it}^{^{\prime }}\overset{d}{\longrightarrow }%
\Sigma _{F}^{-1}\int dB_{\eta }B_{\varepsilon i}^{^{\prime }}=O_{p}\left(
1\right) .
\end{equation*}%
It follows that
\begin{equation*}
\frac{(M_{F^{0}}x_{i})_{t}}{\sqrt{T}}=\frac{x_{it}-\delta _{i}^{\prime
}F_{t}^{0}}{\sqrt{T}}=\frac{x_{it}}{\sqrt{T}}+o_{p}\left( 1\right)
\end{equation*}%
and $\frac{x_{it}}{\sqrt{T}}\overset{d}{\longrightarrow }B_{\varepsilon i}$
as $T\rightarrow \infty $. The limiting distribution of the LS when the
factors are I(0) is the same as when all factors are I(1), except that $%
Q_{i} $ is now asymptotically the same as $B_{\varepsilon i}$. For the FM,
observe that the submatrix $\Omega _{\eta }$ in
\begin{equation*}
\Omega _{bi}=\left[
\begin{array}{cc}
\Omega _{\varepsilon i} & \Omega _{\varepsilon \eta i} \\
\Omega _{\eta \varepsilon i} & \Omega _{\eta }%
\end{array}%
\right]
\end{equation*}%
is a zero matrix since $\eta =\Delta F_{t}^{0}$ is an $I(-1)$ process and
has zero long-run variance. Similarly, $\Omega _{\varepsilon \eta i}$ is
also zero. The submatrix $\Omega _{u\eta i}$ in $\Omega _{u.bi}=\Omega
_{ui}-\Omega _{ubi}\Omega _{bi}^{-1}\Omega _{bui}$ as well as the
submatrices $(%
\begin{array}{ll}
\Delta _{\eta ui} & \Delta _{\eta i}%
\end{array}%
)$ in $(%
\begin{array}{ll}
\Delta _{bui} & \Delta _{bi}%
\end{array}%
)$ are also degenerate because the factors are I(0). Note that $\Omega _{bi}$
is not invertible. Under appropriate choice of bandwidth, see Phillips
(1995), $\Omega _{bi}^{-1}\Omega _{bui}$ can be consistently estimated, so
that FM estimators can be constructed. This argument treats $F_{t}$ as if it
were I(1). If it is known that $F_{t}$ is I(0), we will simply use $F_{t}$
instead of $\Delta F_{t}$ in the FM construction.

\paragraph{ I(1) Regressors, Mixed I(0)/I(1) Factors}

Consider the model
\begin{equation}
y_{it}=x_{it}^{^{\prime }}\beta +\lambda _{1i}^{^{\prime }}F_{1t}+\lambda
_{2i}^{^{\prime }}F_{2t}+u_{it}  \label{mf}
\end{equation}%
where $F_{1t}=\eta _{1t}$ is $r_{1}\times 1$ and $\Delta F_{2t}=\eta _{2t}$
is $r_{2}\times 1$. We again have $M_{F^{0}}x_{i} =x_{i}-F^{0}\delta _{i}$
but $\delta _{i}=\left[
\begin{array}{cc}
\delta _{1i} & \delta _{2i}%
\end{array}
\right]^\prime . $ Then%
\begin{eqnarray*}
\frac{(M_{F^0} x_i)_t}{\sqrt{T}} &=&\frac{x_{it}}{\sqrt{T}}-\frac{1}{\sqrt{T}%
}\left[
\begin{array}{cc}
\delta _{1i}^{^{\prime }} & \delta _{2i}^{^{\prime }}%
\end{array}%
\right] \left[
\begin{array}{c}
F_{1t}^{0} \\
F_{2t}^{0}%
\end{array}%
\right] =\frac{x_{it}}{\sqrt{T}}-\frac{1}{\sqrt{T}}\left( \delta
_{1i}^{^{\prime }}F_{1t}^{0}+\delta _{2i}^{^{\prime }}F_{2t}^{0}\right) \\
&=&\frac{x_{it}}{\sqrt{T}}-\frac{\delta _{2i}^{^{\prime }}F_{2t}^{0}}{\sqrt{T%
}}+o_{p}\left( 1\right)
\end{eqnarray*}%
since $\delta _{1i}=O_{p}\left( 1\right)$, $\delta _{2i}=O_{p}\left(
1\right) $ but $\frac{F^0_{1t}}{\sqrt{T}}=o_p(1)$. The random matrix $Q_i$
involves $B_{\varepsilon i}$ and $B_{2 \eta}$. In the FM correction, the
long run variance $(u_{it}, \Delta x_{it}^{\prime }, \Delta F_{1t}^{\prime
},\Delta F_{2t}^{\prime })^{\prime }$ is degenerate. With an appropriate
choice of bandwidth as in Phillips (1995), the limiting normality still
holds.

\paragraph{Mixed I(1)/I(0) Regressors and I(1) Factors}

Suppose $k_2$ regressors denoted by $x_{2it}$ are I(1), and $k_1$ regressors
denoted by $x_{1it}$ are I(0). Assume $F_{t}$ is I(1) and $u_{it}$ is I(0)
as in (\ref{f1a} ). Consider
\begin{eqnarray*}
y_{it}&=& \alpha_i + x_{1it}^{^{\prime }}\beta _{1}+x_{2it}^{^{\prime
}}\beta _{2}+\lambda _{i}^{\prime }F_{t}+u_{it}  \label{m1} \\
\Delta x_{2it} &=&\varepsilon _{2it}.
\end{eqnarray*}
With the inclusion of an intercept, there is no loss of generality to assume
$x_{1it}$ having a zero mean. For this model, we add the assumption that
\begin{equation}
E(x_{1it} u_{it}) =0  \label{exogenous}
\end{equation}
to rule out simultaneity bias with I(0) regressors. Otherwise $\beta_1$
cannot be consistently estimated. Alternatively, if $u_{it}$ is correlated
with $x_{1it}$, we can project $u_{it}$ onto $x_{1it}$ to obtain the
projection residual and still denote it by $u_{it}$ (with abuse of
notation), and by definition, $u_{it}$ is uncorrelated with $x_{1it}$. But
then $\beta_1$ is no longer the structural parameter. The dynamic least
squares approach by adding $\Delta x_{2it}$ is exactly based on this
argument, with the purpose of more efficient estimation of $\beta_2$.

If one knows which variable is I(0) and which is I(1), the situation is very
simple. The I(1) and I(0) variables are asymptotically orthogonal, we can
separately analyze the distribution of the estimated $\beta_1$ and $\beta_2$.
The estimated $\beta_1$ needs no correction and is asymptotically normal,
and the estimated $\beta_2$ has a distribution as if there is no I(0)
regressors except the intercept. Note that the FM construction for $\hat \beta_2$
is based on the residuals with all regressors included. The rest of analysis
is identical to the situation of all I(1) regressors with an intercept.

In practice, the separation of I(0) or I(1) regressors may not be known in
advance. One can proceed by pretesting to identify the integration order for
each variable, and then apply the above argument. One major purpose of
separating I(0) and I(1) variables is to derive relevant rate of convergence
for the estimated parameters. But if the ultimate purpose is to do
hypothesis testing, there is no need to know the rate of convergence for the
estimator since the scaling factor $n$ or $T$ are cancelled out in the end.
One can proceed as if all regressors are I(1). Then care should be taken
since the long-run covariance matrix is of deficient rank. Phillips (1995)
shows that FM estimators can be constructed with appropriate choice of
bandwidth. Interested readers are referred to Phillips (1995) for details.

Finally, there is the case of mixed I(1)/I(0) regressors and mixed I(1)/I(0)
factors. As explained earlier, I(0) factors do not change the result.
In practice, there is no need to know whether $F^0$ is I(1)
and I(0), since the Cup estimator only depends on $M_{\hat F}$; scaling in $%
\hat F$ does not alter the numerical value of $\hat \beta_{Cup}$.

\section{Monte Carlo Simulations}

In this section, we conduct Monte Carlo experiments to assess the finite
sample properties of the proposed CupBC and CupFM estimators. We also
compare the performance of the proposed estimators with that of LSDV (least
squares dummy variables, i.e., the within group estimator) and 2sFM (2-stage
fully modified which is the CupFM estimator with only one iteration).

Data are generated based on the following design. For $i=1,...,n,$ $%
t=1,...,T,$
\begin{eqnarray*}
y_{it} &=&2x_{it}\text{ }+c\left( \lambda _{i}^{^{\prime }}F_{t}\right)
+u_{it} \\
F_{t} &=&F_{t-1}+\eta _{t} \\
x_{it} &=&x_{it-1}+\varepsilon _{it}
\end{eqnarray*}%
where\footnote{%
Random numbers for error terms, $\left( u_{it},\varepsilon _{it},\eta
_{t}\right) $ are generated by the GAUSS procedure RNDNS. At each
replication, we generate an $nT$ length of random numbers and then split it
into $n$ series so that each series has the same mean and variance.}
\begin{equation}
\left(
\begin{array}{c}
u_{it} \\
\varepsilon _{it} \\
\eta _{t}%
\end{array}%
\right) \overset{iid}{\sim }N\left( \left[
\begin{array}{c}
0 \\
0 \\
0%
\end{array}%
\right] ,\left[
\begin{array}{ccc}
1 & \sigma _{12} & \sigma _{13} \\
\sigma _{21} & 1 & \sigma _{23} \\
\sigma _{31} & \sigma _{32} & 1%
\end{array}%
\right] \right) .  \label{e0}
\end{equation}

We assume a single factor, i.e., $r=1$, $\lambda _{i}$ and $\eta _{t}$ are
generated from i.i.d. $N(\mu _{\lambda },1)$ and $N(\mu _{\eta },1)$
respectively. We set $\mu _{\lambda }=2$ and $\mu _{\eta }=0$. Endogeneity
in the system is controlled by only two parameters, $\sigma _{21}$ and $%
\sigma _{31}.$ The parameter $c$ controls the importance of the global
stochastic trends. We consider $c=\left( 5,10\right),$ $\sigma _{32}=0.4,$ $%
\sigma _{21}=\left( 0,0.2,-0.2\right) $ and $\sigma _{31}=\left(
0,0.8,-0.8\right) .$

The long-run covariance matrix is estimated using the KERNEL procedure in
COINT $2.0$. We use the Bartlett window with the truncation set at five.
Results for other kernels, such as Parzen and quadratic spectral kernels,
are similar and hence not reported. The maximum number of the iteration for
CupBC and CupFM estimators is set to $20$.

Table 1 reports the means and standard deviations (in parentheses) of the
estimators for sample sizes $T=n$ $=\left( 20,40,60,120\right).$ The results
are based on $10,000$ replications. The bias of the LSDV estimator does not
decrease as $\left( n,T\right) $ increases in general. In terms of mean
bias, the CupBC and CupFM are distinctly superior to the LSDV and 2sFM
estimators for all cases considered. The 2sFM estimator is less efficient
than the CupBC and CupFM estimators, as seen by the larger standard
deviations.

To see how the properties of the estimator vary with $n$ and $T$, Table 2
considers $16$ different combinations for $n$ and $T$, each ranging from $20$
to $120$. From Table 2, we see that the LSDV\ and 2sFM\ estimators become
heavily biased when the importance of the common shock is magnified as we
increase $c$ from $5$ to $10$. On the other hand, the CupBC and CupFM
estimators are unaffected by the values of $c$. The results in Table 2 again
indicate that the CupBC and CupFM perform well.

The properties of the $t$-statistic for testing ${\beta =\beta _{0},}$ are
given in Table 3. Here, the LSDV $t$-statistic is the conventional $t$%
-statistic as reported by standard statistical packages. It is clear that
LSDV $t$-statistics and 2sFM $t$-statistics diverge as $\left( n,T\right) $
increases and they are not well approximated by a standard \textrm{N}(0,1)
distribution. The CupBC and CupFM\ $t$-statistics are much better
approximated by a standard \textrm{N}(0,1). Interesting, the performance of
CupBC is no worse than that of CupFM, even though CupBC does the full
modification in the final stage of iteration.

Table 4 shows that, as $n$ and $T$ increases, the biases for the
t-statistics associated with LSDV and 2sFM do not decrease. For CupBC and
CupFM, the biases for the $t$-statistics become smaller (except for a small
number of cases) as $T$ increases for each fixed $n$. As $n$ increases, no
improvement in bias is found. The large standard deviations in the $t$%
-statistics associated with LSDV and 2sFM indicate their poor performance,
especially as $T$ increases. For the CupBC and CupFM, the standard errors
converge to $1.0$ as $n$ and $T$ (especially as $T$) increase.

\section{Conclusion}

This paper develops an asymptotic theory for a panel cointegration
model with unobservable global stochastic trends. Standard least
squares estimator is, in general, inconsistent. In contrast, the
proposed Cup estimator is shown to be consistent (at least
$T$-consistent). In the absence of endogeneity, the Cup estimator is
also $\sqrt{n}T$ consistent. Because we allow the regressors and the
unobservable trends to be endogenous, an asymptotic bias exists for
the Cup estimator.   We further consider two bias-corrected
estimators, CupBC and CupFM, and derive their rate of convergence
and their limiting distributions. We show that these estimators are
$\sqrt{n}T$ consistent and this holds in spite of endogeneity and in
spite of spuriousness induced by unobservable I(1) common shocks. A
simulation study shows that the proposed CupBC and CupFM\ estimators
have good finite sample properties.

\newpage

\appendix

\section*{Appendix}

\renewcommand{\thelem}{A.\arabic{lem}} \renewcommand{\theprop}{A.%
\arabic{prop}} \renewcommand{\thecorollary}{A.\arabic{corollary}} %
\setcounter{lem}{0} \setcounter{prop}{0}

Throughout we use $\left( n,T\right) _{seq}\rightarrow \infty $ to
denote the sequential limit, i.e., $T$ $\rightarrow \infty $ first
and followed by  $n\rightarrow \infty $. We use $MN(0,V)$ to denote
a mixed normal distribution with variance $V$. Let $C$ be the
$\sigma $-field generated by  $\left\{ F_{t}^0\right\} .$ The first
lemma assumes $u_i$ is uncorrelated with $(u_i,F^0)$ for every $i$.
This assumption is relaxed in Lemma \ref{zz}.

\begin{lem}
\label{referee}Suppose that Assumptions 1-5 hold and that $u_i$ is
uncorrelated with $(x_i,F^0)$, then as  $\left( n,T\right)
_{seq}\rightarrow \infty $

\begin{enumerate}
\item[(a)]
\begin{equation*}
\frac{1}{n}\sum_{i=1}^{n}\frac{1}{T^{2}}x_{i}^{^{\prime }}M_{F^0}x_{i}\overset{%
p}{\longrightarrow }\underset{n\rightarrow \infty }{\lim }\frac{1}{n}%
\sum_{i=1}^{n}E\left( \int Q_{i}Q_{i}^{^{\prime }}|C\right) ,
\end{equation*}

\item[(b)]
\begin{equation*}
\frac{1}{\sqrt{n}}\sum_{i=1}^{n}\frac{1}{T}x_{i}^{^{\prime }}M_{F^0}u_{i}%
\overset{d}{\longrightarrow }MN\left( 0,\lim_{n\rightarrow \infty }\frac{1}{n%
}\sum_{i=1}^{n}\Omega _{ui}E\left( \int Q_{i}Q_{i}^{^{\prime }}|C\right)
\right) .
\end{equation*}
\end{enumerate}
\end{lem}

\proof%
Note that
\begin{eqnarray*}
&&\frac{1}{n}\sum_{i=1}^{n}\frac{1}{T^{2}}x_{i}^{^{\prime
}}M_{F^0}x_{i}  =
\frac{1}{n}\sum_{i=1}^{n}\frac{1}{T^{2}}x_{i}^{^{\prime
}}M_{F^0}M_{F^0}x_{i}
= \frac{1}{n}\sum_{i=1}^{n}\frac{1}{T^{2}}\sum_{t=1}^{T}\widetilde{x}_{it}%
\widetilde{x}_{it}^{^{\prime }}
\end{eqnarray*}%
where $\widetilde{x}_{it}=x_{it}-\delta _{i}^{^{\prime }}F^0_{t}$
and \begin{equation*} \delta _{i}=\left( F^{^{0\prime }}F^0\right)
^{-1}F^{^{0\prime }}x_{i}=\left(
\frac{F^{^{0\prime }}F^0}{T^{2}}\right) ^{-1}\frac{1}{T^{2}}%
\sum_{t=1}^{T}F^0_{t}x_{it}^{^{\prime }}\overset{d}{\longrightarrow
}\left( \int B_{\eta }B_{\eta }^{^{\prime }}\right) ^{-1}\int
B_{\eta }B_{\varepsilon i}^{^{\prime }}
\end{equation*}%
see,  e.g., Phillips and Ouliaris (1990).  Thus
\begin{equation*}
\frac{\widetilde{x}_{it}}{\sqrt{T}}=\frac{x_{it}}{\sqrt{T}}-\delta
_{i}^{^{\prime }}\frac{F^0_{t}}{\sqrt{T}}\overset{d}{\longrightarrow }%
B_{\varepsilon i}-\left[ \left( \int B_{\eta }B_{\eta }^{^{\prime
}}\right) ^{-1}\int B_{\eta }B_{\varepsilon i}^{^{\prime }}\right]
^{^{\prime }}B_{\eta }=Q_{i}.
\end{equation*}%
By the continuous mapping theorem
\begin{equation*}
\frac{1}{T^{2}}\sum_{t=1}^{T}\widetilde{x}_{it}\widetilde{x}_{it}^{^{\prime
}}\overset{d}{\longrightarrow }\int Q_{i}Q_{i}^{^{\prime }}=\zeta _{1i}
\end{equation*}%
as $T\rightarrow \infty .$  The variable $\zeta _{1i}$ is
independent across $i$ conditional on $C$, which is an invariant $\sigma $%
-field.\ Thus conditioning on $C$, the law of large numbers for
independent random variables gives,
\begin{equation*}
\frac{1}{n}\sum_{i=1}^{n}\zeta _{1i}\overset{p}{\longrightarrow }%
\lim_{n\rightarrow \infty }\frac{1}{n}\sum_{i=1}^{n}E\left( \zeta
_{1i}|C\right) =\lim_{n\rightarrow \infty
}\frac{1}{n}\sum_{i=1}^{n}E\left( \int Q_{i}Q_{i}^{^{\prime
}}|C\right).
\end{equation*}%
Thus, the sequential limit is
\begin{equation*}
\frac{1}{n}\sum_{i=1}^{n}\frac{1}{T^{2}}x_{i}^{^{\prime }}M_{F^0}x_{i}\overset{%
p}{\longrightarrow }\lim_{n\rightarrow \infty }\frac{1}{n}%
\sum_{i=1}^{n}E\left( \int Q_{i}Q_{i}^{^{\prime }}|C\right)
\end{equation*}%
 This proves part
(a).

Consider (b). Rewrite
\begin{equation*}
\frac{1}{\sqrt{n}}\sum_{i=1}^{n}\frac{1}{T}x_{i}^{^{\prime }}M_{F^0}u_{i}=%
\frac{1}{\sqrt{n}}\sum_{\text{ }i=1}^{n}\frac{1}{T}\sum_{t=1}^{T}\widetilde{x%
}_{it}u_{it}
\end{equation*}%
where $\widetilde{x}_{it}=x_{it}-\delta _{i}^{^{\prime }}F_{t}^{0}$
as before. By assumption, $u_{it}$ is I(0) and is  uncorrelated with
$\widetilde {x}_{it}$. It follows that
\begin{equation*}
\frac{1}{T}\sum_{t=1}^{T}\widetilde{x}_{it}u_{it}\overset{d}{\longrightarrow
}\int Q_{i}dB_{ui}=\xi _{2i}\sim \Omega _{ui}^{1/2}\left( \int
Q_{i}Q_{i}^{^{\prime }}\right) ^{1/2}\times Z
\end{equation*}%
where $Z\sim N\left( 0,I_{k}\right) $ as $T\rightarrow \infty $ for a fixed $%
n$. The variable $\xi _{2i}$ is independent across $i$ conditional
on $C$, which is an invariant $\sigma $-field.\ Thus conditioning on
$C$,
\begin{equation}
\frac{1}{n}\sum_{i=1}^{n}\xi _{2i}\xi _{2i}^{^{\prime }}\overset{p}{%
\longrightarrow }\lim_{n\rightarrow \infty }\frac{1}{n}\sum_{i=1}^{n}E\left(
\xi _{2i}\xi _{2i}^{^{\prime }}|C\right) =\lim_{n\rightarrow \infty }\frac{1%
}{n}\sum_{i=1}^{n}\Omega _{ui}\int E\left( Q_{i}Q_{i}^{^{\prime
}}|C\right). \label{above}
\end{equation}%

Let $I_{i}$ be the $\sigma $ field generated by $\{F^0_{t}\}$ and
$\left( \xi _{21},...,\xi _{2i}\right) $. Then $\left\{ \xi
_{2i},I_{i};i\geq 1\right\} $ is a martingale difference sequence
(MDS) because $\left\{ \xi _{2i}\right\} $ are independent across
$i$ conditional on $C$ and
\begin{equation*}
E\left( \xi _{2i}|I_{i-1}\right) =E\left( \xi _{2i}|C\right) =0.
\end{equation*}%
From $ \sum_{i=1}^{n}\xi _{2i}\xi _{2i}^{^{\prime }}=O_{p}\left(
n\right)$,  the conditional Lindeberg condition in Corollary 3.1 of
Hall and Heyde (1980)
 can be written as%
\begin{equation}
\frac{1}{n}\sum_{i=1}^{n}E\left( \xi _{2i}\xi _{2i}^{^{\prime }}1\left(
\left\Vert \xi _{2i}\right\Vert >\sqrt{n}\delta \right) |I_{i-1}\right)
\overset{p}{\rightarrow }0  \label{lin}
\end{equation}
for all $\delta >0.$  To see (\ref{lin}), notice that
\begin{eqnarray*}
&&\frac{1}{n}\sum_{i=1}^{n}E\left( \xi _{2i}\xi _{2i}^{^{\prime
}}1\left( \left\Vert \xi _{2i}\right\Vert >\sqrt{n}\delta \right)
|I_{i-1}\right) = \frac{1}{n}\sum_{i=1}^{n}E\left( \xi _{2i}\xi
_{2i}^{^{\prime }}1\left( \left\Vert \xi _{2i}\right\Vert
>\sqrt{n}\delta \right) |C\right) .
\end{eqnarray*}%
Without loss of generality we assume $\xi _{2i}$ is a scalar to save
notations. By the Cauchy-Schwarz inequality \[ E\left( \xi
_{2i}^{2}1\left( \left\Vert \xi _{2i}\right\Vert >\sqrt{n}\delta
\right) |C\right) \leq \left\{ E\left( \xi _{2i}^{4}|C\right)
\right\} ^{1/2}\left\{ E\left[ 1\left( \left\Vert \xi _{2i}\right\Vert >%
\sqrt{n}\delta \right) |C\right] \right\} ^{1/2}. \] Furthermore,
\[ E \left[ 1\left( \left\Vert \xi _{2i}\right\Vert >\sqrt{n}\delta \right)|
C \right] \leq \frac{E\left( \xi _{2i}^{2}|C\right) }{n\delta
^{2}}\]  It follows that
\[
\frac{1}{n}\sum_{i=1}^{n}E\left( \xi _{2i}^{2}1\left( \left\Vert \xi
_{i}\right\Vert >\sqrt{n}\delta \right) |C\right)  \leq   \frac 1
{\sqrt{n} \delta}  \left[ \frac{1}{n}\sum_{i=1}^{n}\left [
 E\left( \xi _{2i}^{4}|C\right)
E\left( \xi _{2i}^{2}|C\right) \right ] ^{1/2}  \right]=O_p(
n^{-1/2})
\]
in view of
\begin{equation*}
\frac{1}{n}\sum_{i=1}^{n} [ E\left( \xi _{2i}^{4}|C\right)
E(\xi_{2i}^2|C)] ^{1/2}=O_{p}\left( 1\right).
\end{equation*}
This proves (\ref{lin}). The central limit theorem for martingale
difference sequence, e.g., Corollary 3.1 of Hall and Heyde (1980),
implies that
\begin{equation}
\frac{1}{\sqrt{n}}\sum_{i=1}^{n}\xi _{2i}\overset{d}{\longrightarrow }\left[
\lim_{n\rightarrow \infty }\frac{1}{n}\sum_{i=1}^{n}E\left( \xi _{2i}\xi
_{2i}^{^{\prime }}|C\right) \right] ^{1/2}\times Z
\end{equation}%
where $Z\sim N\left( 0,I\right) $ and $Z$ is independent of $\lim_{n\rightarrow \infty }\frac{1}{n}%
\sum_{i=1}^{n}E\left( \xi _{2i}\xi _{2i}^{^{\prime }}|C\right). $
Note that
\begin{equation*}
\left[ \lim_{n\rightarrow \infty }\frac{1}{n}\sum_{i=1}^{n}E\left( \xi
_{2i}\xi _{2i}^{^{\prime }}|C\right) \right] ^{1/2}=\left(
\lim_{n\rightarrow \infty }\frac{1}{n}\sum_{i=1}^{n}\Omega _{ui}E\left( \int
Q_{i}Q_{i}^{^{\prime }}|C\right) \right) ^{1/2}.
\end{equation*}%
Thus, as $\left( n,T\right) _{seq}\rightarrow \infty ,$ we have
\begin{equation*}
\frac{1}{\sqrt{n}}\frac{1}{T}\sum_{i=1}^{n}\sum_{t=1}^{T}\widetilde{x}%
_{it}u_{it}\overset{d}{\longrightarrow }\left( \lim_{n\rightarrow \infty }%
\frac{1}{n}\sum_{i=1}^{n}\Omega _{ui}E\left( \int
Q_{i}Q_{i}^{^{\prime }}|C\right) \right) ^{1/2}\times Z
\end{equation*}%
which is a mixed normal. The above can be rewritten as
\begin{equation*}
\frac{1}{\sqrt{n}}\frac{1}{T}\sum_{i=1}^{n}\sum_{t=1}^{T}\widetilde{x}%
_{it}u_{it}\overset{d}{\longrightarrow } MN\left( 0,\lim_{n\rightarrow \infty }%
\frac{1}{n}\sum_{i=1}^{n}\Omega _{ui}E\left( \int Q_{i}Q_{i}^{^{\prime
}}|C\right) \right) .
\end{equation*}%
This proves part (b).%
\endproof%

The proofs for Propositions \ref{pro1} and \ref{pro2} (with observable $F$)
follow immediately from Lemma \ref{referee}. Propositions 3 and 4 are proved
in the supplementary appendix of Bai et al.\ (2006).


To derive the limiting distribution for $\widehat{\beta }_{Cup}$, we need
the following lemma. Hereafter, we define $\delta _{nT}=\min \left\{ \sqrt{n}%
,T\right\} .$

\begin{lem}
\label{zz} Suppose Assumptions 1-5 hold. Let $%
Z_{i}=M_{F^{0}}x_{i}-\frac{1}{n}\sum_{k=1}^{n}M_{F^{0}}x_{k}a_{ik}$. Then
as $\left( n,T\right) _{seq}\rightarrow \infty $

\begin{enumerate}
\item[(a)]
\begin{equation*}
\frac{1}{nT^{2}}\sum_{i=1}^{n}Z_{i}^{^{\prime }}Z_{i}\overset{p}{%
\longrightarrow }\text{ }\underset{n\rightarrow \infty }{\lim }\frac{1}{n}%
\sum_{i=1}^{n}E\left( \int R_{ni}R_{ni}^{^{\prime }}|C\right) ,
\end{equation*}

\item[(b)] If $u_{i}$ is uncorrelated with $(x_{i},F^{0})$ for all $i$, then
\begin{equation*}
\frac{1}{\sqrt{n}T}\sum_{i=1}^{n}Z_{i}^{^{\prime }}u_{i}\overset{d}{%
\longrightarrow }MN\left( 0,\text{ }\underset{n\rightarrow \infty }{\lim }%
\frac{1}{n}\sum_{i=1}^{n}\Omega _{ui}E\left( \int R_{ni}R_{ni}^{^{\prime
}}|C\right) \right)
\end{equation*}

\item[(c)] If $u_{i}$ is possibly correlated with $(x_{i},F^{0})$, then
\begin{equation*}
\frac{1}{\sqrt{n}T}\sum_{i=1}^{n}Z_{i}^{^{\prime }}u_{i}- \sqrt{n} \, \theta ^{n} \overset{%
d}{\longrightarrow }MN\left( 0,\text{ }\underset{n\rightarrow \infty }{\lim }%
\frac{1}{n}\sum_{i=1}^{n}\bar{\Omega}_{u.bi}E\left( \int
R_{ni}R_{ni}^{^{\prime }}|C\right) \right)
\end{equation*}%
where%
\begin{eqnarray*}
R_{ni} &=&Q_{i}-\frac{1}{n}\sum_{k=1}^{n}Q_{k}a_{ik}, \\
a_{ik} &=&\lambda _{i}^{^{\prime }}\left( \Lambda ^{^{\prime }}\Lambda
/n\right) ^{-1}\lambda _{k}, \\
Q_{i} &=&B_{\varepsilon i}-\left( \int B_{\varepsilon i}B_{\eta }^{^{\prime
}}\right) \left( \int B_{\eta }B_{\eta }^{^{\prime }}\right) ^{-1}B_{\eta }
\\
\theta ^{n} &=&\frac{1}{n}\sum_{i=1}^{n}\left[ \frac{1}{T}Z_{i}^{^{\prime
}}\left(
\begin{array}{cc}
\Delta \bar{x}_{i} & \Delta F%
\end{array}%
\right) \bar{\Omega}_{bi}^{-1}\bar{\Omega}_{bui}+\left(
\begin{array}{cc}
I_{k} & -\bar{\delta}_{i}^{^{\prime }}%
\end{array}%
\right) \left(
\begin{array}{c}
\bar{\Delta}_{\varepsilon ui}^{+} \\
\bar{\Delta}_{\eta u}^{+}%
\end{array}%
\right) \right]
\end{eqnarray*}%
with $\bar{\delta}_{i}=\left( F^{0^{\prime }}F^{0}\right) ^{-1}F^{0^{\prime
}}\bar{x}_{i},$ and $\bar{x}_{i}=x_{i}-\frac{1}{n}\sum_{k=1}^{n}x_{k}a_{ik}$.
\end{enumerate}
\end{lem}

\noindent \textbf{Proof of (a)}. Recall
\begin{eqnarray*}
M_{F^{0}}x_{i} &=&x_{i}-F^{0}\left( F^{0^{\prime }}F^{0}\right)
^{-1}F^{0^{\prime }}x_{i}=x_{i}-F^{0}\delta _{i}
\end{eqnarray*}%
where
\begin{equation*}
\delta _{i}=\left( F^{0^{\prime }}F^{0}\right) ^{-1}F^{0^{\prime
}}x_{i}=\left( \frac{F^{0^{\prime }}F^{0}}{T^{2}}\right) ^{-1}\frac{1}{T^{2}}%
\sum_{t=1}^{T}F_{t}^{0}x_{it}^{^{\prime }}\overset{d}{\longrightarrow }%
\left( \int B_{\eta }B_{\eta }^{^{\prime }}\right) ^{-1}\int B_{\eta
}B_{\varepsilon i}^{^{\prime }}=\pi _{i}
\end{equation*}%
is a $r\times k$ matrix as $T\rightarrow \infty .$  Write
$\widetilde{x}_{i}=M_{F^0} x_i=x_{i}-F^{0}\delta _{i}$,
a $T\times k$ matrix. Hence%
\begin{eqnarray*}
Z_{i} &=&M_{F^{0}}x_{i}-\frac{1}{n}\sum_{k=1}^{n}M_{F^{0}}x_{k}a_{ik} \\
&=&\left( x_{i}-F^{0}\delta _{i}\right) -\frac{1}{n}\sum_{k=1}^{n}\left(
x_{k}-F^{0}\delta _{k}\right) a_{ik}=\widetilde{x}_{i}-\frac{1}{n}%
\sum_{k=1}^{n}\widetilde{x}_{k}a_{ik}
\end{eqnarray*}%
where $a_{ik}=\lambda _{i}^{^{\prime }}\left( \Lambda ^{^{\prime }}\Lambda
/n\right) ^{-1}\lambda _{k}$ is a scalar and
\begin{equation*}
\frac{\widetilde{x}_{it}}{\sqrt{T}}=\frac{x_{it}}{\sqrt{T}}-\delta
_{i}^{^{\prime }}\frac{F_{t}^{0}}{\sqrt{T}}\overset{d}{\longrightarrow }%
B_{\varepsilon i}-\left[ \left( \int B_{\eta }B_{\eta }^{^{\prime }}\right)
^{-1}\int B_{\eta }B_{\varepsilon i}^{^{\prime }}\right] ^{^{\prime
}}B_{\eta }=Q_{i}
\end{equation*}%
a $k\times 1$ vector, as $%
T\rightarrow \infty .$ It follows that%
\begin{equation*}
\frac{Z_{it}}{\sqrt{T}}\overset{d}{\longrightarrow }Q_{i}-\frac{1}{n}%
\sum_{k=1}^{n}Q_{k}a_{ik}=R_{ni}
\end{equation*}%
and using similar steps in part (a) in Lemma \ref{referee} as $n\rightarrow
\infty ,$
\begin{equation*}
\frac{1}{nT^{2}}\sum_{i=1}^{n}\int R_{ni}R_{ni}^{^{\prime }}\overset{p}{%
\longrightarrow }\text{ }\underset{n\rightarrow \infty }{\lim }\frac{1}{n}%
\sum_{i=1}^{n}E\left( \int R_{ni}R_{ni}^{^{\prime }}|C\right) .
\end{equation*}%
Hence%
\begin{equation*}
\frac{1}{nT^{2}}\sum_{i=1}^{n}Z_{i}^{^{\prime }}Z_{i}\overset{p}{%
\longrightarrow }\text{ }\underset{n\rightarrow \infty }{\lim }\frac{1}{n}%
\sum_{i=1}^{n}E\left( \int R_{ni}R_{ni}^{^{\prime }}|C\right)
\end{equation*}%
as $\left( n,T\right) _{seq}\rightarrow \infty, $ showing (a).

\noindent \textbf{Proof of part (b)}. Notice that
\begin{eqnarray*}
\frac{1}{\sqrt{n}T}\sum_{i=1}^{n}Z_{i}^{^{\prime }}u_{i} &=&\frac{1}{\sqrt{n}%
T}\sum_{i=1}^{n}\left( M_{F^{0}}x_{i}-\frac{1}{n}%
\sum_{k=1}^{n}M_{F^{0}}x_{k}a_{ik}\right) ^{^{\prime }}u_{i} \\
&=&\frac{1}{\sqrt{n}T}\sum_{i=1}^{n}\left( M_{F^{0}}x_{i}\right) ^{^{\prime
}}u_{i}-\frac{1}{\sqrt{n}T}\sum_{i=1}^{n}\left( \frac{1}{n}%
\sum_{k=1}^{n}M_{F^{0}}x_{k}a_{ik}\right) ^{^{\prime }}u_{i} \\
&=&I_{b}+II_{b}.
\end{eqnarray*}%
$I_{b}$ is proved in Lemma \ref{referee}, as $(n,T)_{seq}
\rightarrow \infty,$
\begin{equation*}
I_{b}=\frac{1}{\sqrt{n}T}\sum_{i=1}^{n}\left( M_{F^{0}}x_{i}\right)
^{^{\prime }}u_{i}\overset{d}{\longrightarrow }MN\left( 0,\text{ }\underset{%
n\rightarrow \infty }{\lim }\frac{1}{n}\sum_{i=1}^{n}\Omega _{ui}E\left(
\int Q_{i}Q_{i}^{^{\prime }}|C\right) \right)
\end{equation*}%
if $\widetilde{x}_{it}$ and $u_{it}$ are uncorrelated. Similarly,
for $II_{b}$, we have
\begin{equation*}
\frac{1}{\sqrt{n}T}\sum_{i=1}^{n}\Big(\frac{1}{n}%
\sum_{k=1}^{n}a_{ik}M_{F^{0}}x_{k}\Big)^{^{\prime }}u_{i}\overset{d}{%
\longrightarrow }MN\left( 0,\text{ }\underset{n\rightarrow \infty }{\lim }%
\frac{1}{n}\sum_{i=1}^{n}\Omega _{ui}E\left( C_{ni}|C\right) \right)
\end{equation*}%
where $C_{ni}=\frac{1}{n}\sum_{k=1}^{n}a_{ik}\int Q_{k}Q_{k}^{\prime }$ we
have used the fact that $\frac{1}{n^{2}}\sum_{k=1}^{n}%
\sum_{j=1}^{n}a_{ik}a_{ij}=\frac{1}{n}\sum_{k=1}^{n}a_{ik}$. Thus both $%
I_{b} $ and $II_{b}$ have a proper limiting distribution. These
distributions are dependent since they depend on the same $u_{i}$. We can
also derive their joint limiting distribution. Given the form of $Z_{i}$, it
is easy to show that the above convergences imply part (b).

\textbf{Proof of part (c)}. Now suppose $\widetilde{x}_{it}$\ and $u_{it}$
are correlated. It is known that%
\begin{eqnarray}
\frac{1}{T}\sum_{t=1}^{T}\widetilde{x}_{it}u_{it} &=&\frac{1}{T}%
\sum_{t=1}^{T}\left( x_{it}-\delta _{i}^{^{\prime }}F_{t}^{0}\right) u_{it}=%
\frac{1}{T}\sum_{t=1}^{T}\left(
\begin{array}{cc}
I_{k} & -\delta _{i}^{^{\prime }}%
\end{array}%
\right) \left(
\begin{array}{c}
x_{it} \\
F_{t}^{0}%
\end{array}%
\right) u_{it}  \notag \\
&=&\left(
\begin{array}{cc}
I_{k} & -\delta _{i}^{^{\prime }}%
\end{array}%
\right) \frac{1}{T}\sum_{t=1}^{T}\left(
\begin{array}{c}
x_{it} \\
F_{t}^{0}%
\end{array}%
\right) u_{it}  \notag \\
&&\overset{d}{\longrightarrow }\left(
\begin{array}{cc}
I_{k} & -\pi _{i}^{^{\prime }}%
\end{array}%
\right) \left[ \int \left(
\begin{array}{c}
B_{\varepsilon i}dB_{ui} \\
B_{\eta }dB_{ui}%
\end{array}%
\right) +\left(
\begin{array}{c}
\Delta _{\varepsilon ui} \\
\Delta _{\eta u}%
\end{array}%
\right) \right]  \notag \\
&=&\int Q_{i}dB_{ui}+\left(
\begin{array}{cc}
I_{k} & -\pi _{i}^{^{\prime }}%
\end{array}%
\right) \left(
\begin{array}{c}
\Delta _{\varepsilon ui} \\
\Delta _{\eta u}%
\end{array}%
\right)  \label{xu}
\end{eqnarray}%
as $T\rightarrow \infty $ (e.g., Phillips and Durlauf, 1986). First we note
\begin{eqnarray*}
\int Q_{i}dB_{ui} &=&\int Q_{i}d\left( \Omega _{u.bi}^{1/2}V_{i}+\Omega
_{ubi}\Omega _{bi}^{-1/2}W_{i}\right) \\
&=&\int Q_{i}dB_{u.bi}+\int Q_{i}dB_{bi}^{^{\prime }}\Omega _{bi}^{-1}\Omega
_{bui}
\end{eqnarray*}%
such that
\begin{eqnarray*}
E\left[ \int Q_{i}dV_{i}\right] &=&E\left[ E\left[ \int Q_{i}dV_{i}\right]
|\pi _{i}\right] \\
&=&E\left[ E\left[ \int \left( B_{\varepsilon i}-\pi _{i}^{^{\prime
}}B_{\eta }\right) dV_{i}|\pi _{i}\right] \right] =0.
\end{eqnarray*}%
Note that
\begin{eqnarray*}
&&\frac{1}{T}x_{i}^{^{\prime }}M_{F^{0}}\left(
\begin{array}{cc}
\Delta x_{i} & \Delta F%
\end{array}%
\right) \Omega _{bi}^{-1}\Omega _{bui} \\
&=&\frac{1}{T}\widetilde{x}_{i}^{^{\prime }}\left(
\begin{array}{cc}
\Delta x_{i} & \Delta F%
\end{array}%
\right) \Omega _{bi}^{-1}\Omega _{bui} \\
&=&\left(
\begin{array}{cc}
I_{k} & -\delta _{i}^{^{\prime }}%
\end{array}%
\right) \frac{1}{T}\sum_{t=1}^{T}\left(
\begin{array}{c}
x_{it} \\
F_{t}^{0}%
\end{array}%
\right) \Omega _{bi}^{-1}\Omega _{bui}\left(
\begin{array}{c}
\Delta x_{it} \\
\Delta F_{t}^{0}%
\end{array}%
\right) \\
&&\overset{d}{\longrightarrow }\left(
\begin{array}{cc}
I_{k} & -\pi _{i}^{^{\prime }}%
\end{array}%
\right) \left[ \int \left(
\begin{array}{c}
B_{\varepsilon i} \\
B_{\eta }%
\end{array}%
\right) dB_{bi}^{^{\prime }}\Omega _{bi}^{-1}\Omega _{bui}+\Delta
_{bi}\Omega _{bi}^{-1}\Omega _{bui}\right] .
\end{eqnarray*}

Therefore%
\begin{eqnarray}
&&\frac{1}{T}\widetilde{x}_{i}^{^{\prime }}u_{i}-\left[ \frac{1}{T}%
x_{i}^{^{\prime }}M_{F^{0}}\left(
\begin{array}{cc}
\Delta x_{i} & \Delta F%
\end{array}%
\right) \Omega _{bi}^{-1}\Omega _{bui}+\left(
\begin{array}{cc}
I_{k} & -\delta _{i}^{^{\prime }}%
\end{array}%
\right) \left[ \Delta _{bui}-\Delta _{bi}\Omega _{bi}^{-1}\Omega _{bui}%
\right] \right]  \notag \\
&=&\frac{1}{T}\widetilde{x}_{i}^{^{\prime }}u_{i}-\left[ \frac{1}{T}%
x_{i}^{^{\prime }}M_{F^{0}}\left(
\begin{array}{cc}
\Delta x_{i} & \Delta F%
\end{array}%
\right) \Omega _{bi}^{-1}\Omega _{bui}+\left(
\begin{array}{cc}
I_{k} & -\delta _{i}^{^{\prime }}%
\end{array}%
\right) \Delta _{bui}^{+}\right]  \label{lem2}
\end{eqnarray}%
\begin{equation*}
\overset{d}{\longrightarrow }\Omega _{u.bi}^{1/2}\int Q_{i}dV_{i}\sim \left[
\Omega _{u.bi}^{1/2}\int Q_{i}Q_{i}^{^{\prime }}\right] ^{1/2}\times N\left(
0,I_{k}\right)
\end{equation*}%
where
\begin{equation*}
\Delta _{bui}^{+}=\Delta _{bui}-\Delta _{bi}\Omega _{bi}^{-1}\Omega _{bui}.
\end{equation*}

Let
\begin{equation*}
\theta _{1}^{n}=\frac{1}{n}\sum_{i=1}^{n}\left[ \frac{1}{T}x_{i}^{^{\prime
}}M_{F^{0}}\left(
\begin{array}{cc}
\Delta x_{i} & \Delta F%
\end{array}%
\right) \Omega _{bi}^{-1}\Omega _{bui}+\left(
\begin{array}{cc}
I_{k} & -\delta _{i}^{^{\prime }}%
\end{array}%
\right) \Delta _{bui}^{+}\right] .
\end{equation*}

Then we use similar steps in part (b) in Lemma \ref{referee} to get%
\begin{eqnarray*}
\frac{1}{\sqrt{n}T}\sum_{i=1}^{n}\widetilde{x}_{i}^{^{\prime
}}u_{i}- \sqrt{n} \theta
_{1}^{n} &=&\frac{1}{\sqrt{n}T}\sum_{i=1}^{n}\sum_{t=1}^{T}\widetilde{x}%
_{it}u_{it}-\sqrt{n} \theta _{1}^{n} \\
&&\overset{d}{\longrightarrow }MN\left( 0,\text{ }\underset{n\rightarrow
\infty }{\lim }\frac{1}{n}\sum_{i=1}^{n}\Omega _{u.bi}E\left( \int
Q_{i}Q_{i}^{^{\prime }}|C\right) \right)
\end{eqnarray*}%
as $\left( n,T\right) _{seq}\rightarrow \infty $.

Note $Z_{i}=\widetilde{x}_{i}-\frac{1}{n}\sum_{k=1}^{n}\widetilde{x}%
_{k}a_{ik}$ is a demeaned $\widetilde{x}_{i}$ where $\frac{1}{n}%
\sum_{k=1}^{n}\widetilde{x}_{k}a_{ik}$ is the weighted average of $%
\widetilde{x}_{i}$ with the weight $a_{ik}.$ It follows that
\begin{eqnarray*}
Z_{i} &=&\widetilde{x}_{i}-\frac{1}{n}\sum_{k=1}^{n}\widetilde{x}_{k}a_{ik}
\\
&=&\left( x_{i}-F^{0}\delta _{i}\right) -\frac{1}{n}\sum_{k=1}^{n}\left(
x_{k}-F^{0}\delta _{k}\right) a_{ik} \\
&=&\left( x_{i}-\frac{1}{n}\sum_{k=1}^{n}x_{k}a_{ik}\right) -F^{0}\left(
\delta _{i}-\frac{1}{n}\sum_{k=1}^{n}\delta _{k}a_{ik}\right) ^{^{\prime }}
\\
&=&\bar{x}_{i}-F^{0}\bar{\delta}_{i}^{^{\prime }}
\end{eqnarray*}%
where $\bar{x}_{i}=x_{i}-\frac{1}{n}\sum_{k=1}^{n}x_{k}a_{ik}$ and $\bar{%
\delta}_{i}=\delta _{i}-\frac{1}{n}\sum_{k=1}^{n}\delta _{k}a_{ik}.$

We then can modify (\ref{xu}) as
\begin{eqnarray}
\frac{1}{T}\sum_{t=1}^{T}Z_{it}u_{it} &=&\frac{1}{T}\sum_{t=1}^{T}\left(
\bar{x}_{it}-\bar{\delta }_{i}^{^{\prime }}F_{t}^{0}\right) u_{it}  \notag \\
&=&\frac{1}{T}\sum_{t=1}^{T}\left(
\begin{array}{cc}
I_{k} & -\bar{\delta }_{i}^{^{\prime }}%
\end{array}%
\right) \left(
\begin{array}{c}
\bar{x}_{it} \\
F_{t}^{0}%
\end{array}%
\right) u_{it}=\left(
\begin{array}{cc}
I_{k} & -\bar{\delta }_{i}^{^{\prime }}%
\end{array}%
\right) \frac{1}{T}\sum_{t=1}^{T}\left(
\begin{array}{c}
\bar{x}_{it} \\
F_{t}^{0}%
\end{array}%
\right) u_{it}  \notag \\
&&\overset{d}{\longrightarrow }\left(
\begin{array}{cc}
I_{k} & -\bar{\pi }_{i}^{^{\prime }}%
\end{array}%
\right) \left[  \int \left(
\begin{array}{c}
\bar{B}_{\varepsilon i}dB_{ui} \\
B_{\eta }dB_{ui}%
\end{array}%
\right) +\left(
\begin{array}{c}
\bar{\Delta }_{\varepsilon ui} \\
\bar{\Delta }_{\eta u}%
\end{array}%
\right)   \right]   \notag \\
&=&\int R_{ni}dB_{ui}+\left(
\begin{array}{cc}
I_{k} & -\bar{\pi }_{i}^{^{\prime }}%
\end{array}%
\right) \left(
\begin{array}{c}
\bar{\Delta }_{\varepsilon ui} \\
\bar{\Delta }_{\eta u}%
\end{array}%
\right)  \label{zu}
\end{eqnarray}%
where $\bar{B}_{\varepsilon i}=B_{\varepsilon i}-\frac{1}{n}%
\sum_{k=1}^{n}B_{\varepsilon i}a_{ik}$ and%
\begin{equation*}
\bar{\delta }_{i}=\delta _{i}-\frac{1}{n}\sum_{k=1}^{n}\delta _{k}a_{ik}%
\overset{d}{\longrightarrow }\left( \int B_{\eta }B_{\eta }^{^{\prime
}}\right) ^{-1}\int B_{\eta }\bar{B}_{\varepsilon i}^{^{\prime }}=\bar{\pi }%
_{i}.
\end{equation*}%
The $R_{ni}$ terms appears in the last line in (\ref{zu})  because

\begin{eqnarray*}
\bar{B}_{\varepsilon i}-\bar{\pi}_{i}^{^{\prime }}B_{\eta } &=&\left(
B_{\varepsilon i}-\frac{1}{n}\sum_{k=1}^{n}B_{\varepsilon k}a_{ik}\right)
-\left( \int B_{\eta }B_{\eta }^{^{\prime }}\right) ^{-1}\int B_{\eta
}\left( B_{\varepsilon i}-\frac{1}{n}\sum_{k=1}^{n}B_{\varepsilon
k}a_{ik}\right) ^{^{\prime }}B_{\eta } \\
&=&B_{\varepsilon i}-\left[ \left( \int B_{\eta }B_{\eta }^{^{\prime
}}\right) ^{-1}\int B_{\eta }B_{\varepsilon i}^{^{\prime }}\right] B_{\eta }-%
\frac{1}{n}\sum_{k=1}^{n}\left\{ B_{\varepsilon k}-\left[ \left( \int
B_{\eta }B_{\eta }^{^{\prime }}\right) ^{-1}\int B_{\eta }B_{\varepsilon
k}^{^{\prime }}\right] B_{\eta }\right\} a_{ik} \\
&=&Q_{i}-\frac{1}{n}\sum_{k=1}^{n}Q_{k}a_{ik}=R_{ni}.
\end{eqnarray*}%
Let%
\begin{equation*}
\theta ^{n}=\frac{1}{n}\sum_{i=1}^{n}\left[ \frac{1}{T}Z_{i}^{^{\prime
}}\left(
\begin{array}{cc}
\Delta \bar{x}_{i} & \Delta F%
\end{array}%
\right) \bar{\Omega}_{bi}^{-1}\bar{\Omega}_{bui}+\left(
\begin{array}{cc}
I_{k} & -\bar{\delta}_{i}^{^{\prime }}%
\end{array}%
\right) \bar{\Delta}_{bui}^{+}\right] .
\end{equation*}

Clearly%
\begin{eqnarray*}
\frac{1}{\sqrt{n}T}\sum_{i=1}^{n}Z_{i}^{^{\prime }}u_{i}-\sqrt{n} \, \theta ^{n} &=&%
\frac{1}{\sqrt{n}T}\sum_{i=1}^{n}\left( \widetilde{x}_{i}-\frac{1}{n}%
\sum_{k=1}^{n}\widetilde{x}_{k}a_{ik}\right) ^{^{\prime }}u_{i} \\
&&\overset{d}{\longrightarrow }MN\left( 0,\text{ }\underset{n\rightarrow
\infty }{\lim }\frac{1}{n}\sum_{i=1}^{n}\bar{\Omega}_{u.bi}E\left( \int
R_{ni}R_{ni}^{^{\prime }}|C\right) \right)
\end{eqnarray*}%
as $\left( n,T\rightarrow \infty \right) $ with $R_{ni}=Q_{i}-\frac{1}{n}%
\sum_{k=1}^{n}Q_{k}a_{ik}.$ This proves (c). \endproof

\textbf{Proof of Theorem \ref{theo1}.}

This follows directly from Lemma \ref{zz} as $\left( n,T\right) \rightarrow
\infty $ when $\frac{n}{T}\rightarrow 0$%
\begin{equation*}
\sqrt{n}T\left( \widehat{\beta }_{Cup}-\beta \right) -\sqrt{n}\phi _{nT}%
\overset{d}{\longrightarrow }MN\left( 0,D_{Z}^{-1}\left[ \underset{%
n\rightarrow \infty }{\lim }\frac{1}{n}\sum_{i=1}^{n}\bar{\Omega}%
_{u.bi}E\left( \int R_{ni}R_{ni}^{^{\prime }}|C\right) \right]
D_{Z}^{-1}\right)
\end{equation*}%
where $D_{Z}=\underset{n\rightarrow \infty }{\lim }\frac{1}{n}%
\sum_{i=1}^{n}E\left( \int R_{ni}R_{ni}^{^{\prime }}|C\right) $ and $\phi
_{nT}=\left[ \frac{1}{nT^{2}}\sum_{i=1}^{n}Z_{i}^{^{\prime }}Z_{i}\right]
^{-1}\theta ^{n}.$
\endproof%

\textbf{Proof of Theorem 2 and 3}. The proof for Theorem 2 is similar to that of Theorem 3
below, thus omitted. To prove Theorem 3, we need some preliminary results. First we examine the
limiting distribution of the infeasible FM estimator, $\widetilde{\beta }%
_{CupFM}.$ The endogeneity correction is achieved by modifying the variable $%
y_{it}\,$in (\ref{f1a}) with the transformation%
\begin{equation*}
y_{it}^{+}=y_{it}-\bar{\Omega }_{ubi}\bar{\Omega }_{bi}^{-1}\left(
\begin{array}{c}
\Delta \bar{x}_{it} \\
\Delta F_{t}^{0}%
\end{array}%
\right)
\end{equation*}%
and%
\begin{equation*}
u_{it}^{+}=u_{it}-\bar{\Omega }_{ubi}\bar{\Omega }_{bi}^{-1}\left(
\begin{array}{c}
\Delta \bar{x}_{it} \\
\Delta F_{t}^{0}%
\end{array}%
\right) .
\end{equation*}

By construction $u_{it}^{+}$ has zero long-run covariance with $\left(
\begin{array}{cc}
\Delta \bar{x}_{it}^{^{\prime }} & \Delta F_{t}^{^{0\prime }}%
\end{array}%
\right) ^{^{\prime }}$ and hence the endogeneity can be removed. The serial
correlation correction term has the form%
\begin{eqnarray*}
\bar{\Delta }_{bui}^{+} &=&\left(
\begin{array}{c}
\bar{\Delta }_{\varepsilon ui}^{+} \\
\bar{\Delta }_{\eta u}^{+}%
\end{array}%
\right) =\left(
\begin{array}{ll}
\bar{\Delta }_{bui} & \bar{\Delta }_{bi}%
\end{array}%
\right) \left(
\begin{array}{c}
I_{k} \\
-\bar{\Omega }_{bi}^{-1}\bar{\Omega }_{bui}%
\end{array}%
\right) \\
&=&\bar{\Delta }_{bui}-\bar{\Delta }_{bi}\bar{\Omega }_{bi}^{-1}\bar{\Omega }%
_{bui},
\end{eqnarray*}%
where $\bar{\Delta }_{bui}$ denotes the one-sided long-run covariance
between $u_{it}$ and $\left( \varepsilon _{it},\eta _{t}\right) .$
Therefore, the infeasible FM estimator is%
\begin{equation*}
\widetilde{\beta }_{CupFM}=\left( \sum_{i=1}^{n}x_{i}^{^{\prime
}}M_{F^{0}}x_{i}\right) ^{-1}\sum_{i=1}^{n}\left( x_{i}^{^{\prime
}}M_{F^{0}}y_{i}^{+}-T\left( \bar{\Delta }_{\varepsilon ui}^{+}-\bar{\delta }%
_{i}^{^{\prime }}\bar{\Delta }_{\eta u}^{+}\right) \right)
\end{equation*}%
with $\bar{\delta }_{i}=\left( F^{0^{\prime }}F^{0}\right) ^{-1}F^{0^{\prime
}}\bar{x}_{i}.$

The following Lemma gives the limiting distribution of $\widetilde{\beta }%
_{CupFM}.$

\begin{lem}
\label{fm}Suppose Assumptions in Theorem \ref{theo1} hold. Then as $\left(
n,T\right) _{\func{seq}}\rightarrow \infty $
\begin{equation*}
\sqrt{n}T\left( \widetilde{\beta }_{CupFM}-\beta ^{0}\right) \overset{d}{%
\longrightarrow }MN\left( 0,D_{Z}^{-1}\left[ \underset{n\rightarrow \infty }{%
\lim }\frac{1}{n}\sum_{i=1}^{n}\bar{\Omega}_{u.bi}E\left( \int
R_{ni}R_{ni}^{^{\prime }}|C\right) \right] D_{Z}^{-1}\right) .
\end{equation*}
\end{lem}

\proof%
Let $w_{it}^{+}=\left(
\begin{array}{lll}
u_{it}^{+} & \varepsilon _{it}^{^{\prime }} & \eta ^{^{\prime }}%
\end{array}%
\right) ^{^{\prime }}$ and we have
\begin{equation}
\frac{1}{\sqrt{T}}\sum_{t=1}^{\left[ Tr\right] }w_{it}^{+}\overset{d}{%
\longrightarrow }\left[
\begin{array}{l}
B_{ui}^{+} \\
B_{\varepsilon i} \\
B_{\eta }%
\end{array}%
\right] =\left[
\begin{array}{l}
B_{ui}^{+} \\
B_{bi}%
\end{array}%
\right] =BM\left( \Omega _{i}^{+}\right) \text{ as }T\rightarrow \infty ,
\end{equation}%
where%
\begin{eqnarray*}
B_{bi} &=&\left[
\begin{array}{l}
B_{\varepsilon i} \\
B_{\eta }%
\end{array}%
\right] ,\quad \quad \Omega _{u.bi}=\Omega _{ui}-\Omega _{ubi}\Omega
_{bi}^{-1}\Omega _{bui}, \\
\Omega _{i}^{+} &=&\left[
\begin{array}{ll}
\Omega _{u.bi} & 0 \\
0 & \Omega _{bi}%
\end{array}%
\right] =\left[
\begin{array}{ccc}
\Omega _{u.bi} & 0 & 0 \\
0 & \Omega _{\varepsilon i} & \Omega _{\varepsilon \eta i} \\
0 & \Omega _{\eta \varepsilon i} & \Omega _{\eta }%
\end{array}%
\right] \\
&=&\Sigma ^{+}+\Gamma ^{+}+\Gamma ^{+^{\prime }}, \\
\left[
\begin{array}{l}
B_{ui}^{+} \\
B_{bi}%
\end{array}%
\right] &=&\left[
\begin{array}{cc}
I & -\Omega _{ubi}\Omega _{bi}^{-1} \\
0 & I%
\end{array}%
\right] \left[
\begin{array}{l}
B_{ui} \\
B_{bi}%
\end{array}%
\right] .
\end{eqnarray*}%
Define $\Delta _{i}^{+}=\Sigma _{i}^{+}+\Gamma _{i}^{+}.$ and let $%
u_{1it}^{+}=u_{it}-\Omega _{ubi}\Omega _{bi}^{-1}\left(
\begin{array}{c}
\Delta x_{it} \\
\Delta F_{t}%
\end{array}%
\right) .$ First we notice from (\ref{lem2}) in Lemma \ref{zz} that
\begin{eqnarray}
\zeta _{1iT}^{+} &=&\frac{1}{T}\sum_{t=1}^{T}\widetilde{x}%
_{it}u_{1it}^{+}=\left(
\begin{array}{cc}
I_{k} & -\delta _{i}^{^{\prime }}%
\end{array}%
\right) \frac{1}{T}\sum_{t=1}^{T}\left(
\begin{array}{c}
x_{it} \\
F_{t}^{0}%
\end{array}%
\right) u_{1it}^{+}  \notag \\
&=&\left(
\begin{array}{cc}
I_{k} & -\delta _{i}^{^{\prime }}%
\end{array}%
\right) \left[ \frac{1}{T}\sum_{t=1}^{T}\left(
\begin{array}{c}
x_{it} \\
F_{t}^{0}%
\end{array}%
\right) u_{it}-\frac{1}{T}\sum_{t=1}^{T}\left(
\begin{array}{c}
x_{it} \\
F_{t}^{0}%
\end{array}%
\right) \Omega _{ubi}\Omega _{bi}^{-1}\left(
\begin{array}{c}
\Delta x_{it} \\
\Delta F_{t}^{0}%
\end{array}%
\right) \right]  \notag \\
& \overset{d}{\longrightarrow } &\Omega _{u.bi}^{1/2}\int Q_{i}dV_{i}+\left(
\Delta _{\varepsilon ui}^{+}-\pi _{i}^{^{\prime }}\Delta _{\eta u}^{+}\right)
\label{fm_q}
\end{eqnarray}%
as $T\rightarrow \infty .$ Now let
\begin{equation*}
\zeta _{1iT}^{\ast }=\zeta _{1iT}^{+}-\left( \Delta _{\varepsilon
ui}^{+}-\delta _{i}^{^{\prime }}\Delta _{\eta u}^{+}\right) .
\end{equation*}%
Clearly,
\begin{equation*}
\zeta _{1iT}^{\ast }\overset{d}{\longrightarrow }\Omega _{u.bi}^{1/2}\int
Q_{i}dV_{i}.
\end{equation*}%
Thus,
\begin{eqnarray*}
\frac{1}{\sqrt{n}T}\sum_{i=1}^{n}\left( x_{i}^{^{\prime
}}M_{F^{0}}u_{1i}^{+}-T\left( \Delta _{\varepsilon ui}^{+}-\delta
_{i}^{^{\prime }}\Delta _{\eta u}^{+}\right) \right) &=&\frac{1}{\sqrt{n}T}%
\sum_{i=1}^{n}\left( \sum_{t=1}^{T}\widetilde{x}_{it}u_{1it}^{+}-T\left(
\Delta _{\varepsilon ui}^{+}-\delta _{i}^{^{\prime }}\Delta _{\eta
u}^{+}\right) \right) \\
&&\overset{d}{\longrightarrow }MN\left( 0,\text{ }\underset{n\rightarrow
\infty }{\lim }\frac{1}{n}\sum_{i=1}^{n}\Omega _{u.bi}E\left( \int
Q_{i}Q_{i}^{^{\prime }}|C\right) \right)
\end{eqnarray*}%
as $\left( n,T\right) _{\func{seq}}\rightarrow \infty .$ Next, we modify (%
\ref{fm_q}).
\begin{eqnarray*}
\frac{1}{T}\sum_{t=1}^{T}Z_{it}u_{it}^{+} &=&\frac{1}{T}\sum_{t=1}^{T}\left(
\bar{x}_{it}-\bar{\delta}_{i}^{^{\prime }}F_{t}^{0}\right) u_{it}^{+} \\
&=&\left(
\begin{array}{cc}
I_{k} & -\bar{\delta}_{i}^{^{\prime }}%
\end{array}%
\right) \left[ \frac{1}{T}\sum_{t=1}^{T}\left(
\begin{array}{c}
\bar{x}_{it} \\
F_{t}^{0}%
\end{array}%
\right) u_{it}^{+}-\frac{1}{T}\sum_{t=1}^{T}\left(
\begin{array}{c}
\bar{x}_{it} \\
F_{t}^{0}%
\end{array}%
\right) \Omega _{ubi}\Omega _{bi}^{-1}\left(
\begin{array}{c}
\Delta \bar{x}_{it} \\
\Delta F_{t}^{0}%
\end{array}%
\right) \right] \\
&&\overset{d}{\longrightarrow }\left(
\begin{array}{cc}
I_{k} & -\bar{\pi}_{i}^{^{\prime }}%
\end{array}%
\right) \left\{ \int \left(
\begin{array}{c}
\bar{B}_{\varepsilon i} \\
B_{\eta }%
\end{array}%
\right) dB_{ui}+\left(
\begin{array}{c}
\bar{\Delta}_{\varepsilon ui} \\
\bar{\Delta}_{\eta u}%
\end{array}%
\right) -\left[ \int \left(
\begin{array}{c}
\bar{B}_{\varepsilon i} \\
B_{\eta }%
\end{array}%
\right) dB_{bi}^{^{\prime }}\bar{\Omega}_{bi}^{-1}\bar{\Omega}_{bui}+\bar{%
\Delta}_{bi}\right] \right\} \\
&=&\int R_{ni}dB_{ui}+\left(
\begin{array}{cc}
I_{k} & -\bar{\pi}_{i}^{^{\prime }}%
\end{array}%
\right) \left(
\begin{array}{c}
\bar{\Delta}_{\varepsilon ui} \\
\bar{\Delta}_{\eta u}%
\end{array}%
\right) -\int \left[ R_{ni}dB_{bi}^{^{\prime }}\bar{\Omega}_{bi}^{-1}\bar{%
\Omega}_{bui}+\left(
\begin{array}{cc}
I_{k} & -\bar{\pi}_{i}^{^{\prime }}%
\end{array}%
\right) \left(
\begin{array}{c}
\bar{\Delta}_{\varepsilon i} \\
\bar{\Delta}_{\eta }%
\end{array}%
\right) \bar{\Omega}_{bi}^{-1}\bar{\Omega}_{bui}\right] \\
&=&\bar{\Omega}_{u.bi}^{1/2}\int R_{ni}dV_{i}+\left( \bar{\Delta}%
_{\varepsilon ui}^{+}-\bar{\pi}_{i}^{^{\prime }}\bar{\Delta}_{\eta
u}^{+}\right)
\end{eqnarray*}

Therefore,%
\begin{equation*}
\frac{1}{\sqrt{n}T}\sum_{i=1}^{n}\left( Z_{i}^{^{\prime }}u_{i}^{+}-T\left(
\bar{\Delta}_{\varepsilon ui}^{+}-\bar{\delta}_{i}^{^{\prime }}\bar{\Delta}%
_{\eta u}^{+}\right) \right) \overset{d}{\longrightarrow }MN\left( 0,\text{ }%
\underset{n\rightarrow \infty }{\lim }\frac{1}{n}\sum_{i=1}^{n}\bar{\Omega}%
_{u.bi}E\left( \int R_{ni}R_{ni}^{^{\prime }}|C\right) \right)
\end{equation*}%
as $\left( n,T\right) _{\func{seq}}\rightarrow \infty .$ Then%
\begin{eqnarray*}
&&\sqrt{n}T\left( \widetilde{\beta }_{CupFM}-\beta ^{0}\right) \\
&&\overset{d}{\longrightarrow }MN\left( 0,D_{Z}^{-1}\underset{n\rightarrow
\infty }{\lim }\frac{1}{n}\sum_{i=1}^{n}\Omega _{u.bi}E\left( \int
R_{ni}R_{ni}^{^{\prime }}|C\right) D_{Z}^{-1}\right)
\end{eqnarray*}%
as $\left( n,T\right) _{\func{seq}}\rightarrow \infty $. This proves the
theorem.%
\endproof%

To show $\sqrt{n}T\left( \widehat{\beta }_{CupFM}-\widetilde{\beta }%
_{CupFM}\right) =o_{p}\left( 1\right)$, we need the following lemma.

\begin{lem}
\label{fm1}Under Assumptions of 1-6, we have

\begin{enumerate}
\item[(a)] $\sqrt{n}\left( \widehat{\Delta }_{\varepsilon un}^{+}-\Delta
_{\varepsilon un}^{+}\right) =o_{p}\left( 1\right) , $

\item[(b)] $\frac{1}{\sqrt{n}}\sum_{i=1}^{n}\left( \delta _{i}^{^{\prime }}%
\widehat{\Delta }_{\eta u}^{+}-\delta _{i}^{^{\prime }}\Delta _{\eta
u}^{+}\right) =o_{p}\left( 1\right) , $

\item[(c)] $\frac{1}{\sqrt{n}T}\sum_{i=1}^{n}\left( x_{i}^{^{\prime }}M_{%
\widehat{F}}\widehat{u}_{i}^{+}-x_{i}^{^{\prime }}M_{F^{0}}u_{i}^{+}\right)
=o_{p}\left( 1\right) $

where $\widehat{u}_{it}^{+}=u_{it}-\widehat{\Omega }_{ubi}\widehat{\Omega }%
_{bi}^{-1}\left(
\begin{array}{c}
\Delta x_{it} \\
\Delta \widehat{F}_{t}%
\end{array}%
\right)$ , $\widehat{\Delta }_{\varepsilon un}^{+}=\frac{1}{n}\sum_{i=1}^{n}%
\widehat{\Delta }_{\varepsilon ui}^{+} $ and $\Delta _{\varepsilon un}^{+}=%
\frac{1}{n}\sum_{i=1}^{n}\Delta _{\varepsilon ui}^{+}. $
\end{enumerate}
\end{lem}

Note that the lemma holds when the long run variances are replaced by the
bar versions. Since the proofs are basically the same (as demonstrated in
the proof of Theorem 1), the proof is focused on the variances without the
bar.

\proof
First, note that%
\begin{equation*}
\Delta _{bui}^{+}=\left(
\begin{array}{c}
\Delta _{\varepsilon ui}^{+} \\
\Delta _{\eta u}^{+}%
\end{array}%
\right) =\left(
\begin{array}{ll}
\Delta _{bui} & \Delta _{bi}%
\end{array}%
\right) \left(
\begin{array}{c}
1 \\
-\Omega _{bi}^{-1}\Omega _{bui}%
\end{array}%
\right) =\Delta _{bui}-\Delta _{bi}\Omega _{bi}^{-1}\Omega _{bui}.
\end{equation*}%
Then%
\begin{equation*}
\Delta _{\varepsilon ui}^{+}=\Delta _{\varepsilon ui}-\Delta _{\varepsilon
i}\Omega _{\varepsilon i}^{\ast -1}\Omega _{\varepsilon ui}
\end{equation*}%
where $\Omega _{\varepsilon i}^{\ast -1}$ is the first $k\times k$ block of $%
\Omega _{bi}^{-1}.$ Following the arguments as in the proofs of Theorems 9 and 10
of Hannan (1970) (also see similar result of Moon and Perron (2004)), we
have
\begin{eqnarray*}
&&E\left\Vert \sqrt{n}\left( \widehat{\Delta }_{\varepsilon un}^{+}-\Delta
_{\varepsilon un}^{+}\right) \right\Vert ^{2} \\
&\leq &\underset{i}{\sup }E\left\Vert \widehat{\Delta }_{\varepsilon
ui}^{+}-E\widehat{\Delta }_{\varepsilon ui}^{+}\right\Vert ^{2}+n\,\underset{%
i}{\sup }\left\Vert E\widehat{\Delta }_{\varepsilon ui}^{+}-\Delta
_{\varepsilon ui}^{+}\right\Vert ^{2} \\
&=&O\left( \frac{K}{T}\right) +O\left( \frac{n}{K^{2q}}\right) .
\end{eqnarray*}%
It follows that%
\begin{equation*}
\sqrt{n}\left( \widehat{\Delta }_{\varepsilon un}^{+}-\Delta _{\varepsilon
un}^{+}\right) =O_{p}\left( \max \sqrt{\frac{K}{T}},\sqrt{\frac{n}{K^{2q}}}%
\right) .
\end{equation*}%
From Assumption \ref{band}. $K\backsim n^{b}.$ Then%
\begin{equation*}
\frac{n}{K^{2q}}\backsim \frac{n}{n^{2qb}}=n^{\left( 1-2qb\right)
}\rightarrow 0
\end{equation*}%
if $1<2qb$ or $\frac{1}{2q}<b.$ Next%
\begin{equation*}
\frac{K}{T}\backsim \frac{n^{b}}{T}=\exp \left( \log \left( \frac{n^{b}}{T}%
\right) \right) =\exp \left( b-\frac{\log T}{\log n}\right) \log n=n^{b-%
\frac{\log T}{\log n}}\leq n^{b-\lim \inf \frac{\log T}{\log n}}\rightarrow 0
\end{equation*}%
if $b<\lim \inf \frac{\log T}{\log n}$ by Assumption \ref{band}. Then%
\begin{eqnarray*}
\sqrt{n}\left( \widehat{\Delta }_{\varepsilon un}^{+}-\Delta _{\varepsilon
un}^{+}\right) &=&O_{p}\left( \max \sqrt{\frac{K}{T}},\sqrt{\frac{n}{K^{2q}}}%
\right) \\
&=&o_{p}\left( 1\right)
\end{eqnarray*}%
as required. This proves (a).

To establish (b), we note
\begin{eqnarray*}
\frac{1}{\sqrt{n}}\sum_{i=1}^{n}\left( \delta _{i}^{^{\prime }}\widehat{%
\Delta }_{\eta u}^{+}-\delta _{i}^{^{\prime }}\Delta _{\eta u}^{+}\right)
&=&\left( \frac{1}{n}\sum_{i=1}^{n}\delta _{i}^{^{\prime }}\right) \sqrt{n}%
\left( \widehat{\Delta }_{\eta u}^{+}-\Delta _{\eta u}^{+}\right) \\
&=&O_{p}\left( 1\right) O_{p}\left( \max \left\{ \sqrt{\frac{K}{T}},\sqrt{%
\frac{n}{K^{2q}}}\right\} \right) =o_{p}\left( 1\right)
\end{eqnarray*}%
as required for part (b).

Let $\widetilde{u}_{it}^{+}=u_{it}-\widehat{\Omega }_{ubi}\widehat{\Omega }%
_{bi}^{-1}\left(
\begin{array}{c}
\Delta x_{it} \\
\Delta F_{t}%
\end{array}%
\right) $. Next,
\begin{eqnarray*}
&&\frac{1}{\sqrt{n}T}\sum_{i=1}^{n}\left( x_{i}^{^{\prime }}M_{\widehat{F}}%
\widehat{u}_{i}^{+}-x_{i}^{^{\prime }}M_{F^{0}}u_{i}^{+}\right) \\
&=&\frac{1}{\sqrt{n}T}\sum_{i=1}^{n}\left( x_{i}^{^{\prime }}M_{\widehat{F}}%
\widehat{u}_{i}^{+}-x_{i}^{^{\prime }}M_{\widehat{F}}\widetilde{u}%
_{i}^{+}+x_{i}^{^{\prime }}M_{\widehat{F}}\widetilde{u}_{i}^{+}-x_{i}^{^{%
\prime }}M_{\widehat{F}}u_{i}^{+}+x_{i}^{^{\prime }}M_{\widehat{F}%
}u_{i}^{+}-x_{i}^{^{\prime }}M_{F^{0}}u_{i}^{+}\right) \\
&=&\frac{1}{\sqrt{n}T}\sum_{i=1}^{n}\left( x_{i}^{^{\prime }}M_{\widehat{F}}%
\widetilde{u}_{i}^{+}-x_{i}^{^{\prime }}M_{\widehat{F}}u_{i}^{+}\right) +%
\frac{1}{\sqrt{n}T}\sum_{i=1}^{n}\left( x_{i}^{^{\prime }}M_{\widehat{F}%
}u_{i}^{+}-x_{i}^{^{\prime }}M_{F^{0}}u_{i}^{+}\right) +\frac{1}{\sqrt{n}T}%
\sum_{i=1}^{n}\left( x_{i}^{^{\prime }}M_{\widehat{F}}\widehat{u}%
_{i}^{+}-x_{i}^{^{\prime }}M_{\widehat{F}}\widetilde{u}_{i}^{+}\right) \\
&=&\frac{1}{\sqrt{n}T}\sum_{i=1}^{n}x_{i}^{^{\prime }}M_{\widehat{F}}\left(
\widetilde{u}_{i}^{+}-u_{i}^{+}\right) +\frac{1}{\sqrt{n}T}%
\sum_{i=1}^{n}\left( x_{i}^{^{\prime }}M_{\widehat{F}}-x_{i}^{^{\prime
}}M_{F^{0}}\right) u_{i}^{+}+\frac{1}{\sqrt{n}T}\sum_{i=1}^{n}x_{i}^{^{%
\prime }}M_{\widehat{F}}\left( \widehat{u}_{i}^{+}-\widetilde{u}%
_{i}^{+}\right) \\
&=&\frac{1}{\sqrt{n}T}\sum_{i=1}^{n}x_{i}^{^{\prime }}M_{\widehat{F}}\left(
\widetilde{u}_{i}^{+}-u_{i}^{+}\right) +\frac{1}{\sqrt{n}T}%
\sum_{i=1}^{n}x_{i}^{^{\prime }}\left( M_{\widehat{F}}-M_{F^{0}}\right)
u_{i}^{+}+\frac{1}{\sqrt{n}T}\sum_{i=1}^{n}x_{i}^{^{\prime }}M_{\widehat{F}%
}\left( \widehat{u}_{i}^{+}-\widetilde{u}_{i}^{+}\right) \\
&=&I+II+III.
\end{eqnarray*}%
From the proof of Proposition \ref{known-f} in the supplementary appendix,
\begin{equation*}
II=\frac{1}{\sqrt{n}T}\sum_{i=1}^{n}x_{i}^{^{\prime }}\left( M_{\widehat{F}%
}-M_{F^{0}}\right) u_{i}^{+}=o_{p}\left( 1\right)
\end{equation*}%
if we replace $u_{i}$ by $u_{i}^{+}.$ Let $\Delta b_{i}=\left(
\begin{array}{cc}
\Delta x_{i} & \Delta F%
\end{array}%
\right) $ be a $T\times \left( k+r\right) $ matrix. Consider $I$.%
\begin{eqnarray*}
\frac{1}{\sqrt{n}T}\sum_{i=1}^{n}x_{i}^{^{\prime }}M_{\widehat{F}}\left(
\widetilde{u}_{i}^{+}-u_{i}^{+}\right) &=&\frac{1}{\sqrt{n}T}%
\sum_{i=1}^{n}x_{i}^{^{\prime }}M_{\widehat{F}}\left( u_{i}-\Delta b_{i}%
\widehat{\Omega }_{bi}^{-1}\widehat{\Omega }_{bui}-u_{i}+\Delta b_{i}\Omega
_{bi}^{-1}\Omega _{bui}\right) \\
&=&\frac{1}{\sqrt{n}T}\sum_{i=1}^{n}x_{i}^{^{\prime }}M_{\widehat{F}}\left(
\Delta b_{i}\left( \Omega _{ubi}\Omega _{bi}^{-1}-\widehat{\Omega }_{ubi}%
\widehat{\Omega }_{bi}^{-1}\right) \right) \\
&=&\frac{1}{\sqrt{n}T}\sum_{i=1}^{n}x_{i}^{^{\prime }}\left( I_{T}-\frac{%
\widehat{F}\widehat{F}^{^{\prime }}}{T^{2}}\right) \left( \Delta b_{i}\left(
\Omega _{ubi}\Omega _{bi}^{-1}-\widehat{\Omega }_{ubi}\widehat{\Omega }%
_{bi}^{-1}\right) \right) \\
&=&\frac{1}{\sqrt{n}T}\sum_{i=1}^{n}x_{i}^{^{\prime }}\Delta b_{i}\left(
\Omega _{ubi}\Omega _{bi}^{-1}-\widehat{\Omega }_{ubi}\widehat{\Omega }%
_{bi}^{-1}\right) \\
&&-\frac{1}{\sqrt{n}T}\sum_{i=1}^{n}x_{i}^{^{\prime }}\frac{\widehat{F}%
\widehat{F}^{^{\prime }}}{T^{2}}\left( \Delta b_{i}\left( \Omega
_{ubi}\Omega _{bi}^{-1}-\widehat{\Omega }_{ubi}\widehat{\Omega }%
_{bi}^{-1}\right) \right) \\
&=&I_{c}+II_{c}.
\end{eqnarray*}

Along the same lines as the proofs of Theorems 9 and 10 of Hannan (1970), we
can show that
\begin{eqnarray*}
\underset{i}{\sup }E\left\Vert \widehat{\Omega }_{ubi}\widehat{\Omega }%
_{bi}^{-1}-\Omega _{ubi}\Omega _{bi}^{-1}\right\Vert ^{2} &=&O\left( \frac{K%
}{T}\right) +O\left( \frac{1}{K^{2q}}\right) .
\end{eqnarray*}%
Then we have%
\begin{eqnarray*}
\Omega _{ubi}\Omega _{bi}^{-1}-\widehat{\Omega }_{ubi}\widehat{\Omega }%
_{bi}^{-1} =O_{p}\left( Max\left\{ \sqrt{\frac{K}{T}},\sqrt{\frac{1}{K^{2q}}}%
\right\} \right) .
\end{eqnarray*}%
and%
\begin{eqnarray*}
\frac{1}{\sqrt{n}}\sum_{i=1}^{n}\left\Vert \Omega _{ubi}\Omega _{bi}^{-1}-
\widehat{\Omega }_{ubi}\widehat{\Omega }_{bi}^{-1}\right\Vert ^{2} &=&\sqrt{n%
}\frac{1}{n}\sum_{i=1}^{n}\left\Vert \Omega _{ubi}\Omega _{bi}^{-1}-\widehat{%
\Omega }_{ubi}\widehat{\Omega }_{bi}^{-1}\right\Vert ^{2} \\
&\leq &\sqrt{n}\underset{i}{\sup }\left\Vert \Omega _{ubi}\Omega _{bi}^{-1}-%
\widehat{\Omega }_{ubi}\widehat{\Omega }_{bi}^{-1}\right\Vert ^{2} \\
&=&\sqrt{n}\left[ O_{p}\left( Max\left\{ \sqrt{\frac{K}{T}},\sqrt{\frac{1}{%
K^{2q}}}\right\} \right) \right] ^{2}.
\end{eqnarray*}

For $I_{c.}$, by the Cauchy Schwarz inequality, \
\begin{eqnarray*}
&&\left\Vert I_{c}\right\Vert =\left\Vert \frac{1}{\sqrt{n}T}%
\sum_{i=1}^{n}x_{i}^{^{\prime }}\Delta b_{i}\left( \Omega _{ubi}\Omega
_{bi}^{-1}-\widehat{\Omega }_{ubi}\widehat{\Omega }_{bi}^{-1}\right)
\right\Vert \\
&\leq &\left( \sqrt{n}\frac{1}{n}\sum_{i=1}^{n}\left\Vert \frac{%
x_{i}^{^{\prime }}\Delta b_{i}}{T}\right\Vert ^{2}\right) ^{1/2}\left( \frac{%
1}{\sqrt{n}}\sum_{i=1}^{n}\left\Vert \Omega _{ubi}\Omega _{bi}^{-1}-\widehat{%
\Omega }_{ubi}\widehat{\Omega }_{bi}^{-1}\right\Vert ^{2}\right) ^{1/2} \\
&\leq &\left[ O_{p}\left( \sqrt{n}\right) \right] ^{1/2}\left( \sqrt{n}%
\right) ^{1/2}O_{p}\left( Max\left\{ \sqrt{\frac{K}{T}},\sqrt{\frac{1}{K^{2q}%
}}\right\} \right) \\
&=&O_{p}\left( \sqrt{n}\right) O_{p}\left( Max\left\{ \sqrt{\frac{K}{T}},%
\sqrt{\frac{1}{K^{2q}}}\right\} \right)
\end{eqnarray*}

Similarly,%
\begin{eqnarray*}
\left\Vert II_{c.}\right\Vert &=&\left\Vert \frac{1}{\sqrt{n}T}%
\sum_{i=1}^{n}x_{i}^{^{\prime }}\frac{\widehat{F}\widehat{F}^{^{\prime }}}{%
T^{2}}\left( \Delta b_{i}\left( \Omega _{ubi}\Omega _{bi}^{-1}-\widehat{%
\Omega }_{ubi}\widehat{\Omega }_{bi}^{-1}\right) \right) \right\Vert \\
&=&\left\Vert \frac{1}{\sqrt{n}}\sum_{i=1}^{n}\frac{x_{i}^{^{\prime }}%
\widehat{F}}{T^{2}}\frac{\widehat{F}^{^{\prime }}\Delta b_{i}}{T}\left(
\Omega _{ubi}\Omega _{bi}^{-1}-\widehat{\Omega }_{ubi}\widehat{\Omega }%
_{bi}^{-1}\right) \right\Vert \\
&\leq &\left\Vert \left( \sqrt{n}\frac{1}{n}\sum_{i=1}^{n}\left\Vert \frac{%
x_{i}^{^{\prime }}\widehat{F}}{T^{2}}\frac{\widehat{F}^{^{\prime }}\Delta
b_{i}}{T}\right\Vert ^{2}\right) ^{1/2}\left( \frac{1}{\sqrt{n}}%
\sum_{i=1}^{n}\left\Vert \Omega _{ubi}\Omega _{bi}^{-1}-\widehat{\Omega }%
_{ubi}\widehat{\Omega }_{bi}^{-1}\right\Vert ^{2}\right) ^{1/2}\right\Vert \\
&=&O_{p}\left( \sqrt{n}\right) O_{p}\left( Max\left\{ \sqrt{\frac{K}{T}},%
\sqrt{\frac{1}{K^{2q}}}\right\} \right)
\end{eqnarray*}%
Combining $I_{c.}$ and $II_{c},$ we have%
\begin{eqnarray*}
\frac{1}{\sqrt{n}T}\sum_{i=1}^{n}x_{i}^{^{\prime }}M_{\widehat{F}}\left(
\widehat{u}_{i}^{+}-\widetilde{u}_{i}^{+}\right) &=&O_{p}\left( \sqrt{n}%
\right) O_{p}\left( Max\left\{ \sqrt{\frac{K}{T}},\sqrt{\frac{1}{K^{2q}}}%
\right\} \right) \\
&=&O_{p}\left( Max\left\{ \sqrt{\frac{nK}{T}},\sqrt{\frac{n}{K^{2q}}}%
\right\} \right)
\end{eqnarray*}%
Recall $K\backsim n^{b}$ and $\lim \inf \frac{\log T}{\log n}>1$ from
Assumption \ref{band}. It follows that, as in Moon and Perron (2004)%
\begin{eqnarray*}
\frac{nK}{T} &\backsim &\frac{n^{b+1}}{T}=\exp \left( \log \left( \frac{%
n^{b+1}}{T}\right) \right) =\exp \left( b+1-\frac{\log T}{\log n}\right)
\log n \\
&=&n^{b+1-\frac{\log T}{\log n}}\leq n^{b+1-\lim \inf \frac{\log T}{\log n}%
}\rightarrow 0
\end{eqnarray*}%
by Assumption \ref{band} and $b<\lim \inf \frac{\log T}{\log n}-1.$ Also note%
\begin{equation*}
\frac{n}{K^{2q}}\backsim \frac{n}{n^{2qb}}=n^{\left( 1-2qb\right)
}\rightarrow 0
\end{equation*}%
by Assumption \ref{band} and $\frac{1}{2q}<b$. Therefore%
\begin{equation*}
\frac{1}{\sqrt{n}T}\sum_{i=1}^{n}x_{i}^{^{\prime }}M_{\widehat{F}}\left(
\widehat{u}_{i}^{+}-\widetilde{u}_{i}^{+}\right) =O_{p}\left( Max\left\{
\sqrt{\frac{nK}{T}},\sqrt{\frac{n}{K^{2q}}}\right\} \right) =o_{p}\left(
1\right) .
\end{equation*}%
Let
\begin{equation*}
\Delta \widehat{b}_{i}=\left(
\begin{array}{cc}
\Delta x_{i} & \Delta \widehat{F}%
\end{array}%
\right) .
\end{equation*}%
Note that
\begin{equation*}
\Delta b_{i}-\Delta \widehat{b}_{i}=\left(
\begin{array}{cc}
\Delta x_{i} & \Delta F%
\end{array}%
\right) -\left(
\begin{array}{cc}
\Delta x_{i} & \Delta \widehat{F}%
\end{array}%
\right) =\left(
\begin{array}{cc}
0 & \Delta F-\Delta \widehat{F}%
\end{array}%
\right) .
\end{equation*}

Consider III.
\begin{eqnarray*}
\frac{1}{\sqrt{n}T}\sum_{i=1}^{n}x_{i}^{^{\prime }}M_{\widehat{F}}\left(
\widehat{u}_{i}^{+}-\widetilde{u}_{i}^{+}\right) &=&\frac{1}{\sqrt{n}T}%
\sum_{i=1}^{n}x_{i}^{^{\prime }}M_{\widehat{F}}\left( u_{i}-\Delta \widehat{b%
}_{i}\widehat{\Omega }_{bi}^{-1}\widehat{\Omega }_{bui}-u_{i}+\Delta b_{i}%
\widehat{\Omega }_{bi}^{-1}\widehat{\Omega }_{bui}\right) \\
&=&\frac{1}{\sqrt{n}T}\sum_{i=1}^{n}x_{i}^{^{\prime }}M_{\widehat{F}}\left(
\Delta b_{i}-\Delta \widehat{b}_{i}\right) \widehat{\Omega }_{bi}^{-1}%
\widehat{\Omega }_{bui} \\
&=&\frac{1}{\sqrt{n}T}\sum_{i=1}^{n}x_{i}^{^{\prime }}M_{\widehat{F}}\left(
\Delta F-\Delta \widehat{F}\right) \widehat{\Omega }_{bi}^{-1}\widehat{%
\Omega }_{bui}.
\end{eqnarray*}%
We use Lemma 12.3 in Bai (2005) to get
\begin{equation*}
\frac{1}{nT}\sum_{i=1}^{n}x_{i}^{^{\prime }}M_{\widehat{F}}\left( \Delta
F-\Delta \widehat{F}\right) =O_{p}\left( \widehat{\beta }-\beta ^{0}\right)
+O_{p}\left( \frac{1}{\min \left( n,T\right) }\right) .
\end{equation*}%
It follows that
\begin{eqnarray*}
\frac{1}{\sqrt{n}T}\sum_{i=1}^{n}x_{i}^{^{\prime }}M_{\widehat{F}}\left(
\Delta F-\Delta \widehat{F}\right) &=&\sqrt{n}\left[ O_{p}\left( \widehat{%
\beta }-\beta ^{0}\right) +O_{p}\left( \frac{1}{\min \left( n,T\right) }%
\right) \right] \\
&=&\sqrt{n}O_{p}\left( \frac{1}{T}\right) +O_{p}\left( \frac{\sqrt{n}}{\min
\left( n,T\right) }\right)=o_p(1)
\end{eqnarray*}%
since $\frac{n}{T}\rightarrow 0$ as $(n,T)\rightarrow \infty .$ Collecting $%
I-III$ we prove (c).%
\endproof%

\begin{prop}
\label{infe-fe} Under Assumptions \ref{fl}-\ref{band},
\begin{equation*}
\sqrt{n}T\left( \widehat{\beta }_{CupFM}-\widetilde{\beta }_{CupFM}\right)
=o_{p}\left( 1\right) .
\end{equation*}
\end{prop}

\proof%
To save the notations, we only show that results with $x_i$ in place of $%
\bar{x}_{i}$ and $\delta_i$ in place of of $\bar{\delta}_{i}$ since the
steps are basically the same. In the supplementary appendix, it is shown
that (see the proof of Proposition \ref{known-f})
\begin{equation*}
\left( \frac{1}{nT^{2}}\sum_{i=1}^{n}x_{i}^{^{\prime }}M_{\widehat{F}%
}x_{i}\right) =\left( \frac{1}{nT^{2}}\sum_{i=1}^{n}x_{i}^{^{\prime
}}M_{F^{0}}x_{i}\right) +o_{p}\left( 1\right) .
\end{equation*}%
Then%
\begin{eqnarray*}
&&\sqrt{n}T\left( \widehat{\beta }_{CupFM}-\widetilde{\beta }_{CupFM}\right)
\\
&=&\left( \frac{1}{nT^{2}}\sum_{i=1}^{n}x_{i}^{^{\prime
}}M_{F^{0}}x_{i}\right) ^{-1}\frac{1}{\sqrt{n}T}\left\{
\begin{array}{c}
\sum_{i=1}^{n}\left( x_{i}^{^{\prime }}M_{\widehat{F}}\widehat{u}%
_{i}^{+}-T\left( \widehat{\Delta }_{\varepsilon ui}^{+}-\delta
_{i}^{^{\prime }}\widehat{\Delta }_{\eta u}^{+}\right) \right) \\
-\sum_{i=1}^{n}\left( x_{i}^{^{\prime }}M_{F^{0}}u_{i}^{+}-T\left( \Delta
_{\varepsilon ui}^{+}-\delta _{i}^{^{\prime }}\Delta _{\eta u}^{+}\right)
\right)%
\end{array}%
\right\} +o_{p}\left( 1\right) \\
&=&\left( \frac{1}{nT^{2}}\sum_{i=1}^{n}x_{i}^{^{\prime
}}M_{F^{0}}x_{i}\right) ^{-1}\frac{1}{\sqrt{n}T}\left\{
\begin{array}{c}
\sum_{i=1}^{n}\left( x_{i}^{^{\prime }}M_{\widehat{F}}\widehat{u}%
_{i}^{+}-x_{i}^{^{\prime }}M_{F^{0}}u_{i}^{+}\right) \\
-nT\left( \widehat{\Delta }_{\varepsilon un}^{+}-\Delta _{\varepsilon
un}^{+}\right) -T\sum_{i=1}^{n}\left( \delta _{i}^{^{\prime }}\widehat{%
\Delta }_{\eta u}^{+}-\delta _{i}^{^{\prime }}\Delta _{\eta u}^{+}\right)%
\end{array}%
\right\} +o_{p}\left( 1\right) \\
&=&\left( \frac{1}{nT^{2}}\sum_{i=1}^{n}x_{i}^{^{\prime
}}M_{F^{0}}x_{i}\right) ^{-1}\left\{
\begin{array}{c}
\frac{1}{\sqrt{n}T}\sum_{i=1}^{n}\left( x_{i}^{^{\prime }}M_{\widehat{F}}%
\widehat{u}_{i}^{+}-x_{i}^{^{\prime }}M_{F^{0}}u_{i}^{+}\right) \\
-\sqrt{n}\left( \widehat{\Delta }_{\varepsilon un}^{+}-\Delta _{\varepsilon
un}^{+}\right) -\frac{1}{\sqrt{n}}\sum_{i=1}^{n}\left( \delta _{i}^{^{\prime
}}\widehat{\Delta }_{\eta u}^{+}-\delta _{i}^{^{\prime }}\Delta _{\eta
u}^{+}\right)%
\end{array}%
\right\} +o_{p}\left( 1\right)
\end{eqnarray*}%
where $\widehat{\Delta }_{\varepsilon un}^{+}=\frac{1}{n}\sum_{i=1}^{n}%
\widehat{\Delta }_{\varepsilon ui}^{+}$ and$\Delta _{\varepsilon un}^{+}=%
\frac{1}{n}\sum_{i=1}^{n}\Delta _{\varepsilon ui}^{+}.$ Finally using Lemma %
\ref{fm1},
\begin{equation*}
\sqrt{n}T\left( \widehat{\beta }_{CupFM}-\widetilde{\beta }_{CupFM}\right)
=o_{p}\left( 1\right) .
\end{equation*}%
\endproof%

\paragraph{Proof of Theorem \protect\ref{fesi-fm}:}

This follows directly from Proposition \ref{infe-fe}. \endproof

\paragraph{Proof of Proposition \protect\ref{Fhat}:}

In the supplementary appendix, it is shown that
\begin{equation*}
\frac 1 T \sum_{t=1}^T \| \hat F_t -H F_t^0 \|^2 = T \,O_p(\|\hat \beta
-\beta^0\|^2) +O_p(\frac 1 n) + O_p( \frac 1 {T^2}).
\end{equation*}
From $\sqrt{n} T(\hat \beta -\beta^0)=O_p(1 )$, the first term on the right
hand side is $O_p(1/(nT))$, which is dominated by $O(1/n)$. \endproof

\newpage \oddsidemargin -0.7in

\begin{landscape}
\begin{table}[tbp] \centering%

\begin{tabular}{ccccc|cccc|cccc}
\multicolumn{13}{c}{\textbf{Table 1: Mean Bias and Standard Deviation of
Estimators}} \\ \hline
& \multicolumn{4}{c}{$\sigma _{31}=0$} & \multicolumn{4}{c}{$\sigma _{31}=0.8
$} & \multicolumn{4}{c}{$\sigma _{31}=-0.8$} \\
& {\small LSDV} & 2s{\small FM} & CupBC & Cup{\small FM} & {\small LSDV} & 2s%
{\small FM} & CupBC & Cup{\small FM} & {\small LSDV} & 2sFM & CupBC & Cup%
{\small FM} \\ \hline
$\sigma _{21}=0$ &  &  &  &  &  &  &  &  &  &  &  &  \\
n,T=20 & {\small 1.352} & {\small 0.349} & {\small 0.030} & {\small 0.030} & -%
{\small 0.712} & {\small 0.257} & {\small 0.000} & {\small 0.000} & {\small %
2.216} & {\small -0.086} & {\small 0.030} & {\small 0.030} \\
& {\small (1.559)} & {\small (0.387)} & {\small (0.030)} & {\small (0.029)}
& {\small (1.505)} & {\small (0.372)} & {\small (0.030)} & {\small (0.029)}
& {\small (1.524)} & {\small (0.394)} & {\small (0.029)} & {\small (0.029)}
\\
n,T=40 & {\small 3.371} & -{\small 0.719} & {\small -0.000} & {\small -0.000}
& {\small 2.761} & {\small -0.246} & {\small -0.000} & {\small -0.000} &
{\small 1.010} & {\small -0.371} & {\small -0.000} & {\small -0.000} \\
& {\small (1.139)} & {\small (0.225)} & {\small (0.009)} & {\small (0.009)}
& {\small (1.529)} & {\small (0.227)} & {\small (0.010)} & {\small (0.009)}
& {\small (1.124)} & {\small (0.217)} & {\small (0.009)} & {\small (0.009)}
\\
n,T=60 & {\small -2.006} & {\small 0.094} & {\small -0.000} & {\small -0.000}
& -{\small 1.393} & {\small 0.038} & {\small -0.000} & {\small -0.000} &
{\small -1.073} & {\small 0.199} & -{\small 0.000} & {\small -0.000} \\
& {\small (0.920)} & {\small (0.138)} & {\small (0.005)} & {\small (0.005)}
& {\small (0.915)} & {\small (0.139)} & {\small (0.005)} & {\small (0.005)}
& {\small (0.929)} & {\small (0.138)} & {\small (0.005)} & {\small (0.005)}
\\
n,T=120 & {\small 0.204} & -{\small 0.064} & {\small -0.000} & {\small -0.000}
& {\small 0.548} & {\small -0.062} & {\small -0.020} & {\small 0.015} &
{\small -0.163} & {\small -0.061} & {\small 0.018} & {\small -0.000} \\
& {\small (0.645)} & {\small (0.056)} & {\small (0.018)} & {\small (0.002)}
& {\small (0.646)} & {\small (0.056)} & {\small (0.002)} & {\small (0.002)}
& {\small (0.643)} & {\small (0.056)} & {\small (0.002)} & {\small (0.002)}
\\
$\sigma _{21}=0.2$ &  &  &  &  &  &  &  &  &  &  &  &  \\
n,T=20 & {\small 4.333} & {\small 0.317} & {\small -0.119} & {\small 0.332} &
{\small 2.258} & {\small 0.129} & {\small -0.158} & {\small 0.293} & {\small %
4.903} & {\small -0.220} & {\small -0.117} & {\small 0.322} \\
& {\small (1.584)} & {\small (0.385)} & {\small (0.030)} & {\small (0.029)}
& {\small (1.529)} & {\small (0.382)} & {\small (0.031)} & {\small (0.029)}
& {\small (1.614)} & {\small (0.396)} & {\small (0.030)} & {\small (0.028)}
\\
n,T=40 & {\small 4.567} & {\small -0.768} & {\small -0.113} & {\small 0.100} &
{\small 4.051} & -{\small 0.333} & {\small -0.117} & {\small 0.101} &
{\small 1.964} & {\small -0.376} & {\small -0.115} & {\small 0.102} \\
& {\small (1.133)} & {\small (0.223)} & {\small (0.010)} & {\small (0.009)}
& {\small (1.153)} & {\small (0.227)} & {\small (0.010)} & {\small (0.009)}
& {\small (1.120)} & {\small (0.216)} & {\small (0.010)} & {\small (0.009)}
\\
n,T=60 & {\small -1.100} & {\small 0.109} & {\small -0.071} & {\small 0.045} &
{\small -0.337} & {\small 0.082} & {\small -0.067} & {\small 0.049} &
{\small 0.032} & {\small 0.150} & -{\small 0.065} & {\small 0.051} \\
& {\small (0.923)} & {\small (0.138)} & {\small (0.005)} & {\small (0.005)}
& {\small (0.925)} & {\small (0.139)} & {\small (0.005)} & {\small (0.005)}
& {\small (0.938)} & {\small (0.140)} & {\small (0.005)} & {\small (0.005)}
\\
n,T=120 & {\small 0.696} & -{\small 0.059} & {\small 0.000} & {\small 0.178} &
{\small 1.161} & {\small -0.070} & {\small -0.017} & {\small 0.017} &
{\small 0.151} & {\small -0.026} & {\small 0.017} & -{\small 0.017} \\
& {\small (0.648)} & {\small (0.055)} & {\small 0.018} & {\small (0.002)} &
{\small (0.649)} & {\small (0.055)} & {\small (0.002)} & {\small (0.002)} &
{\small (0.646)} & {\small (0.055)} & {\small (0.002)} & {\small (0.002)} \\
$\sigma _{21}=-0.2$ &  &  &  &  &  &  &  &  &  &  &  &  \\
n,T=20 & {\small -1.600} & {\small 0.376} & {\small 0.179} & {\small -0.274} &
{\small -3.763} & {\small 0.331} & {\small 0.151} & {\small -0.291} & -%
{\small 0.754} & {\small -0.049} & {\small 0.169} & {\small -0.274} \\
& {\small (1.588)} & {\small (0.393)} & {\small (0.031)} & {\small (0.029)}
& {\small (1.593)} & {\small (0.345)} & {\small (0.031)} & {\small (0.029)}
& {\small (1.603)} & {\small (0.394)} & {\small (0.031)} & {\small (0.029)}
\\
n,T=40 & {\small 2.086} & {\small -0.653} & {\small 0.105} & {\small -0.108} &
{\small 0.812} & -{\small 0.077} & {\small 0.101} & {\small -0.113} &
{\small -0.353} & {\small -0.313} & {\small 0.096} & {\small -0.112} \\
& {\small (1.144)} & {\small (0.225)} & {\small (0.010)} & {\small (0.009)}
& {\small (1.141)} & {\small (0.223)} & {\small (0.010)} & {\small (0.009)}
& {\small (1.128)} & {\small (0.218)} & {\small (0.010)} & {\small (0.009)}
\\
n,T=60 & {\small -2.850} & {\small 0.008} & {\small 0.055} & {\small -0.062} &
{\small -2.178} & {\small -0.018} & {\small 0.058} & {\small -0.058} &
{\small -1.872} & {\small 0.236} & {\small 0.056} & {\small -0.060} \\
& {\small (0.917)} & {\small (0.142)} & {\small (0.005)} & {\small (0.005)}
& {\small (0.905)} & {\small (0.136)} & {\small (0.005)} & {\small (0.005)}
& {\small (0.921)} & {\small (0.138)} & {\small (0.005)} & {\small (0.005)}
\\
n,T=120 & -{\small 0.501} & {\small 0.000} & {\small 0.000} & {\small 0.000} &
-{\small 0.175} & -{\small 0.000} & {\small -0.018} & {\small 0.017} &
{\small -0.839} & {\small 0.029} & {\small 0.000} & {\small -0.000} \\
& {\small (0.650)} & {\small (0.057)} & {\small (0.002)} & {\small (0.018)}
& {\small (0.646)} & {\small (0.057)} & {\small (0.002)} & {\small (0.002)}
& {\small (0.654)} & {\small (0.058)} & {\small (0.002)} & {\small (0.002)}
\\ \hline
\multicolumn{13}{l}{Note: (a) The Mean biases here have been multiplied by
100.} \\
\multicolumn{13}{l}{(b) $c=5,$ $\sigma _{32}=0.4.$}%
\end{tabular}%


\end{table}%

\end{landscape}

\begin{table}[tbp] \centering%
$%
\begin{tabular}{ccccc|cccc}
\multicolumn{9}{c}{\textbf{Table 2: Mean Bias and Standard Deviation }} \\
\multicolumn{9}{c}{\textbf{of Estimators for Different n and T}} \\ \hline
&  & $c=5$ &  &  &  & $c=10$ &  &  \\ \hline
(n,T) & LSDV & 2sFM & CupBC & CupFM & LSDV & 2sFM & CupBC & CupFM \\ \hline
$(20,20)$ & 2.258 & 0.129 & -0.158 & 0.293 & 1.538 & 0.275 & -0.158 & 0.294
\\
& (1.594) & (0.382) & (0.031) & (0.028) & (3.186) & (0.771) & (0.031) &
(0.029) \\
$(20,40)$ & 4.832 & -0.426 & -0.067 & 0.107 & 8.141 & -0.006 & -0.067 & 0.106
\\
& (1.692) & (0.288) & (0.014) & (0.014) & (3.186) & (0.566) & (0.014) &
(0.014) \\
$(20,60)$ & 0.460 & 0.282 & -0.019 & -0.058 & -0.105 & 0.0561 & -0.186 &
0.058 \\
& (1.560) & (0.206) & (0.009) & (0.009) & (3.121) & (0.412) & (0.009) &
(0.009) \\
$(20,120)$ & 3.018 & 0.040 & 0.010 & 0.021 & -6.550 & 0.067 & 0.010 & 0.021
\\
& (1.572) & (0.123) & (0.005 & (0.005) & (3.144) & (0.245) & (0.005) &
(0.004) \\
$(40,20)$ & 4.012 & -0.566 & -0.225 & 0.320 & 5.092 & -1.087 & -0.226 & 0.320
\\
& (1.126) & (0.280) & (0.0218) & (0.019) & (2.252) & (0.593) & (0.021) &
(0.019) \\
$(40,40)$ & 4.051 & -0.332 & -0.117 & 0.101 & 6.616 & -0.622 & -0.117 & 0.101
\\
& (1.153) & (0.227) & (0.010) & (0.009) & (2.305) & (0.454) & (0.010) &
(0.009) \\
$(40,60)$ & 1.818 & 0.114 & -0.055 & 0.051 & 2.628 & 0.248 & -0.055 & 0.051
\\
& (1.098) & (0.158) & (0.007) & (0.006) & (2.196) & (0.317) & (0.007) &
(0.006) \\
$(40,120)$ & 1.905 & -0.090 & -0.010 & 0.015 & 3.303 & -0.178 & -0.010 &
0.015 \\
& (1.111) & (0.087) & (0.003) & (0.003) & (2.243) & (0.187) & (0.003) &
(0.003) \\
$(60,20)$ & 3.934 & -0.317 & -0.294 & 0.295 & 4.989 & -0.544 & -0.294 & 0.295
\\
& (0.921) & (0.249) & (0.018) & (0.017) & (1.841) & (0.497) & (0.014) &
(0.016) \\
$(60,40)$ & 2.023 & 0.110 & -0.125 & 0.108 & 2.573 & 0.267 & -0.125 & 0.109
\\
& (0.923) & (0.187) & (0.009) & (0.008) & (1.296) & (0.027) & (0.009) &
(0.008) \\
$(60,60)$ & -0.337 & 0.082 & -0.067 & 0.049 & -1.666 & 0.191 & -0.067 & 0.049
\\
& (0.925) & (0.139) & (0.005) & (0.005) & (1.850) & (0.279) & (0.005) &
(0.005) \\
$(60,120)$ & -1.168 & 0.109 & -0.015 & 0.015 & -2.839 & -0.223 & -0.014 &
0.015 \\
& (0.923) & (0.075) & (0.003) & (0.003) & (1.847) & (0.151) & (0.003) &
(0.003) \\
$(120,20)$ & 2.548 & -0.151 & -0.304 & 0.294 & 2.236 & -0.203 & -0.304 &
0.294 \\
& (0.651) & (0.182) & (0.014) & (0.011) & (1.303) & (0.362) & (0.014) &
(0.011) \\
$(120,40)$ & 1.579 & -0.026 & -0.013 & 0.001 & 1.678 & 0.000 & -0.133 & 0.112
\\
& (0.661) & (0.137) & (0.006) & (0.005) & (1.321) & (0.279) & (0.006) &
(0.005) \\
$(120,60)$ & 0.764 & 0.004 & -0.077 & 0.013 & 0.539 & 0.061 & -0.077 & 0.048
\\
& (0.634) & (0.100) & (0.004) & (0.004) & (1.267) & (0.199) & (0.004) &
(0.004) \\
$(120,120)$ & 1.161 & -0.070 & -0.017 & 0.017 & 1.823 & -0.134 & -0.017 &
0.018 \\
& (0.649) & (0.055) & (0.002) & (0.002) & (1.298) & (0.111) & (0.002) &
(0.002) \\ \hline\hline
\multicolumn{9}{l}{(a) The Mean biases here have been multiplied by 100.} \\
\multicolumn{9}{l}{(b) $\sigma _{21}=0.2,$ $\sigma _{31}=0.8,$ and $\sigma
_{32}=0.4.$}%
\end{tabular}%
$%
\end{table}%

\bigskip
\begin{landscape}

\begin{table}[tbp] \centering%
$%
\begin{tabular}{ccccc|cccc|cccc}
\multicolumn{13}{c}{\textbf{Table 3: Mean Bias and Standard Deviation of
t-statistics}} \\ \hline
& \multicolumn{2}{c}{} & $\sigma _{31}=0$ &  & \multicolumn{4}{c}{$\sigma
_{31}=0.8$} & \multicolumn{4}{c}{$\sigma _{31}=-0.8$} \\
& {\small LSDV} & 2s{\small FM} & CupBC & Cup{\small FM} & {\small LSDV} & 2s%
{\small FM} & CupBC & Cup{\small FM} & {\small LSDV} & 2s{\small FM} & CupBC
& Cup{\small FM} \\ \hline
$\sigma _{21}=0$ &  &  &  &  &  &  &  &  &  &  &  &  \\
n,T=20 & {\small 0.036} & {\small 0.006} & 0.016 & 0.016 & {\small 0.006} &
{\small 0.0224} & {\small 0.001} & {\small 0.001} & {\small 0.041} & {\small %
-0.001} & {\small 0.019} & {\small 0.019} \\
& {\small (2.414)} & {\small (2.445)} & {\small (1.531)} & {\small (1.502)}
& {\small (2.527)} & {\small (2.449)} & {\small (1.529)} & {\small (1.503)}
& {\small (2.534)} & {\small (2.455)} & {\small (1.515)} & {\small (1.491)}
\\
n,T=40 & {\small 0.092} & {\small -0.036} & {\small -0.007} & {\small -0.006}
& {\small 0.074} & {\small -0.052} & {\small -0.012} & {\small -0.011} &
{\small 0.019} & {\small 0.008} & {\small -0.006} & {\small -0.005} \\
& {\small (3.576)} & {\small (2.589)} & {\small (1.276)} & {\small (1.256)}
& {\small (3.592)} & {\small (2.618)} & {\small (1.273)} & {\small (1.254)}
& {\small (3.588)} & {\small (2.581)} & {\small (1.278)} & {\small (1.217)}
\\
n,T=60 & {\small -0.098} & {\small 0.016} & {\small -0.019} & {\small -0.019}
& {\small -0.036} & {\small -0.016} & {\small -0.011} & {\small -0.011} &
{\small -0.060} & {\small 0.045} & {\small -0.009} & {\small -0.009} \\
& {\small (4.346)} & {\small (2.647)} & {\small (1.182)} & {\small (1.169)}
& {\small (4.325)} & {\small (2.640)} & {\small (1.189)} & {\small (1.178)}
& {\small (4.315)} & {\small (2.644)} & {\small (1.182)} & {\small (1.169)}
\\
n,T=120 & {\small 0.046} & -{\small 0.019} & {\small -0.003} & {\small -0.003%
} & {\small 0.099} & -{\small 0.019} & {\small -0.075} & {\small 0.102} &
{\small -0.088} & -{\small 0.040} & {\small 0.068} & {\small -0.011} \\
& {\small (6.093)} & {\small (2.696)} & {\small (1.101)} & {\small (1.096)}
& {\small (6.089)} & {\small (2.661)} & {\small (1.118)} & {\small (1.094)}
& {\small (6.095)} & {\small (2.705)} & {\small (1.120)} & {\small (1.095)}
\\
$\sigma _{21}=0.2$ &  &  &  &  &  &  &  &  &  &  &  &  \\
n,T=20 & {\small 0.104} & {\small 0.040} & {\small 0.001} & {\small 0.185} &
{\small 0.070} & {\small 0.037} & {\small -0.013} & {\small 0.188} & {\small %
0.105} & {\small 0.033} & {\small 0.004} & {\small 0.181} \\
& {\small (2.508)} & {\small (2.454)} & {\small (1.558)} & {\small (1.497)}
& {\small (2.529)} & {\small (2.453)} & {\small (1.561)} & {\small (1.442)}
& {\small (2.539)} & {\small (2.465)} & {\small (1.543)} & {\small (1.483)}
\\
n,T=40 & {\small 0.149} & -{\small 0.013} & {\small -0.081} & {\small 0.140}
& {\small 0.134} & {\small -0.022} & {\small -0.085} & {\small 0.142} &
{\small 0.059} & -{\small 0.003} & {\small -0.081} & {\small 0.143} \\
& {\small (3.563)} & {\small (2.597)} & {\small (1.304)} & {\small (1.252)}
& {\small (3.578)} & {\small (2.639)} & {\small (1.307)} & {\small (1.252)}
& {\small (3.578)} & {\small (2.612)} & {\small (1.314)} & {\small (1.258)}
\\
n,T=60 & {\small -0.032} & {\small 0.039} & {\small -0.100} & {\small 0.115}
& {\small 0.027} & {\small 0.013} & {\small -0.094} & {\small 0.123} &
{\small 0.011} & {\small 0.038} & {\small -0.087} & {\small 0.127} \\
& {\small (4.357)} & {\small (2.651)} & {\small (1.209)} & {\small (1.167)}
& {\small (4.357)} & {\small (2.647)} & {\small (1.215)} & {\small (1.174)}
& {\small (4.325)} & {\small (2.646)} & {\small (1.204)} & {\small (1.162)}
\\
n,T=120 & {\small 0.049} & -{\small 0.016} & {\small 0.003} & {\small 0.002}
& {\small 0.097} & {\small -0.019} & {\small -0.059} & {\small 0.114} &
{\small 0.012} & -{\small 0.029} & {\small 0.062} & -{\small 0.109} \\
& {\small (6.060)} & {\small (2.640)} & {\small (1.096)} & {\small (1.092)}
& {\small (6.084)} & {\small (2.645)} & {\small (1.115)} & {\small (1.093)}
& {\small (6.043)} & {\small (2.635)} & {\small (1.111)} & {\small (1.089)}
\\
$\sigma _{21}=-0.2$ &  &  &  &  &  &  &  &  &  &  &  &  \\
n,T=20 & {\small -0.031} & {\small -0.013} & {\small 0.029} & {\small -0.155}
& {\small -0.064} & {\small 0.005} & {\small 0.125} & {\small -0.166} &
{\small -0.031} & {\small -0.029} & {\small 0.027} & {\small -0.152} \\
& {\small (2.519)} & {\small (2.456)} & {\small (1.559)} & {\small (1.497)}
& {\small (2.528)} & {\small (2.439)} & {\small (1.556)} & {\small (1.498)}
& {\small (2.538)} & {\small (2.458)} & {\small (1.556)} & {\small (1.498)}
\\
n,T=40 & {\small 0.033} & {\small -0.068} & {\small 0.067} & {\small -0.153}
& {\small -0.005} & {\small -0.071} & {\small 0.061} & {\small -0.162} &
{\small -0.035} & {\small -0.021} & {\small 0.058} & {\small -0.159} \\
& {\small (3.586)} & {\small (2.593)} & {\small (1.312)} & {\small (1.255)}
& {\small (3.597)} & {\small (2.618)} & {\small (1.305)} & {\small (1.248)}
& {\small (3.588)} & {\small (2.574)} & {\small (1.305)} & {\small (1.252)}
\\
n,T=60 & {\small -0.162} & {\small 0.002} & {\small 0.062} & {\small -0.154}
& {\small -0.093} & {\small -0.035} & {\small 0.067} & {\small -0.146} &
{\small -0.114} & {\small 0.028} & {\small 0.067} & {\small -0.147} \\
& {\small (4.335)} & {\small (2.657)} & {\small (1.212)} & {\small (1.169)}
& {\small (4.283)} & {\small (2.633)} & {\small (1.210)} & {\small (1.168)}
& {\small (4.308)} & {\small (2.643)} & {\small (1.206)} & {\small (1.166)}
\\
n,T=120 & {\small -0.066} & {\small 0.001} & {\small 0.007} & {\small 0.007}
& {\small -0.010} & {\small 0.022} & {\small -0.062} & {\small 0.117} &
{\small -0.111} & -{\small 0.004} & {\small 0.077} & {\small -0.104} \\
& {\small (6.098)} & {\small (2.679)} & {\small (1.106)} & {\small (1.106)}
& {\small (6.152)} & {\small (2.577)} & {\small (1.116)} & {\small (1.092)}
& {\small (6.119)} & {\small (2.691)} & {\small (1.125)} & {\small (1.101)}
\\ \hline
\multicolumn{13}{l}{Note: (a) $c=5$, $\sigma _{32}=0.4.$} \\
\multicolumn{13}{l}{}%
\end{tabular}%
$%
\end{table}%

\bigskip \end{landscape}

\begin{table}[tbp] \centering%
$%
\begin{tabular}{ccccc|cccc}
\multicolumn{9}{c}{\textbf{Table 4: Mean Bias and Standard Deviation }} \\
\multicolumn{9}{c}{\textbf{of t-statistics for Different n and T}} \\ \hline
&  & $c=5$ &  &  &  & $c=10$ &  &  \\ \hline
(n,T) & LSDV & 2sFM & CupBC & CupFM & LSDV & 2sFM & CupBC & CupFM \\ \hline
$(20,20)$ & 0.070 & 0.037 & -0.013 & 0.169 & 0.036 & 0.030 & -0.013 & 0.169
\\
& (2.529) & (2.453) & (1.561) & (1.497) & (2.532) & (2.562) & (1.560) &
(1.496) \\
$(20,40)$ & 0.130 & -0.007 & -0.009 & 0.110 & 0.106 & -0.011 & -0.009 & 0.110
\\
& (3.539) & (1.863) & (1.313) & (1.286) & (3.541) & (1.896) & (1.313) &
(1.286) \\
$(20,60)$ & 0.029 & 0.009 & 0.015 & 0.085 & 0.009 & 0.003 & 0.016 & 0.085 \\
& (4.303) & (1.553) & (1.253) & (1.239) & (4.305) & (1.569) & (1.253) &
(1.239) \\
$(20,120)$ & -0.090 & 0.015 & 0.057 & 0.064 & -0.105 & 0.013 & 0.057 & 0.064
\\
& (6.131) & (1.222) & (1.156) & (1.151) & (6.132) & (1.220) & (1.156) &
(1.151) \\
$(40,20)$ & 0.119 & -0.015 & -0.086 & 0.242 & 0.073 & -0.019 & -0.086 & 0.241
\\
& (2.518) & (3.376) & (1.549) & (1.443) & (2.520) & (3.610) & (1.549) &
(1.443) \\
$(40,40)$ & 0.134 & -0.022 & -0.085 & 0.142 & 0.100 & -0.026 & -0.085 & 0.142
\\
& (3.578) & (2.639) & (1.307) & (1.252) & (3.580) & (2.739) & (1.307) &
(1.252) \\
$(40,60)$ & 0.113 & 0.012 & -0.048 & 0.109 & 0.085 & 0.008 & -0.047 & 0.109
\\
& (4.328) & (2.164) & (1.209) & (1.177) & (4.329) & (2.222) & (1.209) &
(1.176) \\
$(40,120)$ & 0.133 & -0.014 & -0.007 & 0.059 & 0.113 & -0.019 & -0.007 &
0.059 \\
& (6.097) & (1.519) & (1.131) & (1.123) & (6.098) & (1.535) & (1.131) &
(1.123) \\
$(60,20)$ & 0.123 & 0.005 & -0.161 & 0.276 & 0.067 & -0.002 & -0.160 & 0.276
\\
& (2.521) & (4.042) & (1.579) & (1.424) & (2.524) & (4.409) & (1.579) &
(1.425) \\
$(60,40)$ & 0.100 & 0.069 & -0.109 & 0.192 & 0.059 & 0.065 & -0.109 & 0.192
\\
& (3.532) & (3.206) & (1.352) & (1.272) & (3.534) & (3.375) & (1.352) &
(1.272) \\
$(60,60)$ & 0.027 & 0.013 & -0.094 & 0.123 & -0.006 & 0.010 & -0.094 & 0.122
\\
& (4.426) & (2.613) & (1.215) & (1.174) & (4.359) & (2.751) & (1.215) &
(1.174) \\
$(60,120)$ & -0.020 & 0.031 & -0.024 & 0.077 & -0.044 & 0.030 & -0.025 &
0.077 \\
& (6.131) & (1.866) & (1.118) & (1.104) & (6.132) & (1.902) & (1.118) &
(1.104) \\
$(120,20)$ & 0.139 & 0.044 & -0.243 & 0.386 & 0.060 & 0.063 & -0.243 & 0.386
\\
& (2.478) & (5.269) & (1.681) & (1.404) & (2.479) & (5.969) & (1.681) &
(1.404) \\
$(120,40)$ & 0.135 & 0.037 & -0.186 & 0.268 & 0.078 & 0.040 & -0.186 & 0.268
\\
& (3.588) & (4.369) & (1.366) & (1.233) & (3.589) & (4.706) & (1.366) &
(1.233) \\
$(120,60)$ & 0.099 & 0.011 & -0.162 & 0.174 & 0.052 & 0.004 & -0.162 & 0.174
\\
& (4.272) & (3.683) & (1.249) & (1.166) & (4.273) & (3.902) & (1.249) &
(1.167) \\
$(120,120)$ & 0.097 & -0.189 & -0.589 & 0.114 & 0.063 & -0.027 & -0.059 &
0.114 \\
& (6.084) & (2.645) & (1.115) & (1.093) & (6.086) & (2.741) & (1.115) &
(1.093) \\ \hline\hline
\multicolumn{9}{l}{(a) $\sigma _{21}=0.2,$ $\sigma _{31}=0.8,$ and $\sigma
_{32}=0.4.$}%
\end{tabular}%
$%
\end{table}%


\end{document}